\documentclass[twoside]{article}
\usepackage{latexsym}
\usepackage{amsmath}
\usepackage{amssymb}
\usepackage[dvips]{graphicx}

\newcommand{\real}{{\mathbb R}}
\newcommand{\vp}{\varphi}
\newcommand{\vpz}{\varphi^z}
\newcommand{\vpxi}{\varphi^{xy}}
\newcommand{\eps}{\varepsilon}
\newcommand{\tzer}{e^x}
\newcommand{\tone}{\tau_1}
\newcommand{\ttwo}{\tau_2}
\newcommand{\nzer}{e^y}
\newcommand{\none}{\nu_1}
\newcommand{\ntwo}{\nu_2}

\newcommand{\ol}{\overline}
\newcommand{\td}{\tilde}

\newcommand{\dive}{{\rm div}}

\newcommand{\qed}{\hfill $\Box$}
\newcommand{\df}{\partial}
\newcommand{\dfi}{\df_{\nu_i}}
\newcommand{\dfo}{\df_y}
\newcommand{\dfii}{\df^2_{\nu_i}}
\newcommand{\dft}{\df_{\tau_i}}
\newcommand{\dsns}{\df^2_{\nu_i t_1}}
\newcommand{\dsnt}{\df^2_{\nu_i t_2}}
\newcommand{\dsst}{\df^2_{t_1t_2}}
\newcommand{\hess}{\nabla^2_{\nu_i t_1 t_2}}
\newcommand{\vhess}{\nabla^2_{y\, t_1 t_2}}
\newcommand{\vvhess}{\nabla^2_{\nu_2 t_1 t_2}}
\newcommand{\nablav}{\nabla_{\!\nu_i t_1 t_2}}
\newcommand{\nablaz}{\nabla_{\! y\, t_1 t_2}}
\newcommand{\vs}{\vspace{.3cm}}

\newcommand{\bstep}[2]{\vspace{.5cm}\noindent {\sc Step #1.-- } { 
#2}\vspace{.5cm}} 

\newcommand{\aplimsup}{\mathop{\rm ap\,lim\,sup}}
\newcommand{\apliminf}{\mathop{\rm ap\,lim\,inf}}

\makeatletter
\@addtoreset{equation}{section}
\makeatother

\newtheorem{thm}{Theorem}[section]
\newtheorem{lm}[thm]{Lemma}

\begin{document}

\null
\vspace{1cm}
\begin{center}
{\Large \bf LOCAL CALIBRATIONS FOR MINIMIZERS OF THE \\

\

MUMFORD-SHAH FUNCTIONAL \\

\

WITH A TRIPLE JUNCTION}

\vspace{1.5cm}

{\large Maria Giovanna \sc Mora}

\vspace{.6cm}

{\large S.I.S.S.A.}\\ \vspace{.1cm}
{\large via Beirut 2-4, 34014 Trieste, Italy}\\ \vspace{.1cm}
{\large e-mail: \tt mora@sissa.it}

\vspace{3cm}

\begin{minipage}[t]{11.5cm}
{\normalsize
\noindent
{\bf Abstract.} We prove that, if $u$ is a function satisfying all Euler
conditions for the Mumford-Shah functional and the discontinuity set of $u$ is
given by three line segments meeting at the origin with equal angles, then
there exists a neighbourhood $U$ of the origin such that $u$ is a minimizer of
the Mumford-Shah functional on $U$ with respect to its own boundary conditions
on $\partial U$. The proof is obtained by using the calibration
method.}\vspace{.3cm}

 \vspace{.3cm}

\end{minipage}
\end{center}

\setcounter{page}{0}
\thispagestyle{empty}

\vfill

\begin{center}
{\large Ref. S.I.S.S.A. 39/2001/M  (May 2001)}
\end{center}

\pagebreak
\clearpage
\phantom{a}
\setcounter{page}{0}
\thispagestyle{empty}

\newpage
\clearpage
\pagestyle{myheadings}

\null
\vspace{.5cm}
\begin{center}
{\large \bf LOCAL CALIBRATIONS FOR MINIMIZERS OF THE MUMFORD-SHAH 

\vspace{.5cm}

FUNCTIONAL WITH A TRIPLE JUNCTION}
\vspace{1cm}

{\normalsize  Maria Giovanna \sc Mora}

\vspace{1cm}

\begin{minipage}[t]{11.5cm}
{\footnotesize
\noindent
{\bf Abstract.} We prove that, if $u$ is a function satisfying all Euler
conditions for the Mumford-Shah functional and the discontinuity set of $u$ is
given by three line segments meeting at the origin with equal angles, then
there exists a neighbourhood $U$ of the origin such that $u$ is a minimizer of
the Mumford-Shah functional on $U$ with respect to its own boundary conditions
on $\partial U$. The proof is obtained by using the calibration method.}
\end{minipage}
\end{center}

\vspace{1cm}

\section{Introduction}

The Mumford-Shah functional was proposed in \cite{Mum-Sha2} to approach image
segmentation problems and it can be written, in the ``homogeneous'' version in
dimension two, as
\begin{equation}\label{ms}
\int_{\Omega}|\nabla u(x,y)|^2dx\, dy +{\cal H}^1(S_u),
\end{equation}
where $\Omega$ is a bounded open subset of $\real^2$ with a Lipschitz
boundary, ${\cal H}^1$ is the one-dimensional Hausdorff measure, $u$ is the
unknown function in the space $SBV(\Omega)$ of special functions of bounded
variation in $\Omega$, $S_u$ is the set of essential discontinuity points of
$u$, while $\nabla u$ denotes its approximate gradient (see
\cite{Amb-Fus-Pal}).

This paper deals with local minimizers of (\ref{ms}) with given
boundary values. More precisely, we say that $u$ is a {\it Dirichlet
minimizer\/} of (\ref{ms}) in $\Omega$ if $u$ belongs to $SBV(\Omega)$ and
satisfies the inequality
$$\int_{\Omega}|\nabla u(x,y)|^2dx\, dy +{\cal H}^1(S_u)\leq
\int_{\Omega}|\nabla v(x,y)|^2dx\, dy +{\cal H}^1(S_v)$$
for every $v\in SBV(\Omega)$ with the same trace as $u$ on $\df\Omega$.

Considering different classes of infinitesimal variations, one can show that,
if $u$ is a Dirichlet minimizer of (\ref{ms}) in $\Omega$, then the following equilibrium
conditions (which can be globally called the Euler conditions for (\ref{ms}))
are satisfied (see \cite{Amb-Fus-Pal} and \cite{Mum-Sha2}):  
\begin{itemize}
	\item $u$ is harmonic on $\Omega\setminus S_u$;
	\item the normal derivative of $u$ vanishes on both sides of $S_u$,
where $S_u$ is a regular curve;
	\item the curvature of $S_u$ (where defined) is equal to the
difference of the squares of the tangential derivatives of $u$ on both sides
of $S_u$;
	\item if $S_u$ is locally the union of finitely many
regular arcs, then $S_u$ can present only two kinds of singularities: 
either a regular arc ending at some point, the so-called ``crack-tip'', or
three regular arcs meeting with equal angles of $2\pi/3$, the so-called
``triple junction''.  
\end{itemize} 
However, since the functional (\ref{ms}) is not
convex, the Euler conditions are not sufficient for the Dirichlet minimality
of $u$.

In \cite{Mor-Mor} it has been proved that, if $S_u$ is an analytic curve
connecting two points of $\df\Omega$ (hence, $S_u$ has no singular points),
then the Euler conditions are also sufficient for the
Dirichlet minimality in small domains. In this paper we prove that, if $S_u$
is given by three line segments meeting at the origin with equal angles of $2\pi/3$
(i.e., $S_u$ is a rectilinear triple junction), the same conclusion holds; in
other words, for every $(x,y)\in \Omega$, there is an open neighbourhood $U$ of
$(x,y)$ such that $u$ is a Dirichlet minimizer of (\ref{ms}) in $U$. Since for
$(x,y)\neq (0,0)$ this fact follows from the result in \cite{Mor-Mor}, the
interesting case is when we restrict the functional to a neighbourhood of the
triple point $(0,0)$.

The precise statement of the result is the following.

\begin{thm}\label{tj}
Let $\Omega:=B(0,1)$ be the open disc in $\real^2$ with radius $1$ centred at
the origin, and let $(A_0,A_1,A_2)$ be the partition of $\Omega$ defined as
follows: 
$$A_i:=\left\{ (r\cos\theta,r\sin\theta)\in\Omega: \ 0\leq r<1, \
\frac{2}{3}\pi(2-i)<\theta\leq  \frac{2}{3}\pi(3-i)\right\} \qquad \forall \,
i=0,1,2.$$ Let $S_{i,j}:=\ol{A_i}\cap\ol{A_j}$ for every $i<j$.
Let $u_i\in C^2(\ol{A_i})$ be a harmonic function in $A_i$, satisfying the
Neumann conditions on $\df A_i\cap\Omega$ and such that $|\nabla u_i|=|\nabla
u_j|$ on $S_{i,j}$ for every $i<j$. 
If $u$ is the function in $SBV(\Omega)$
defined by $u:= u_i$ a.e. in each $A_i$ and $u_0(0,0)<u_1(0,0)<u_2(0,0)$, then
there exists a neighbourhood $U$ of the origin such that $u$ is a Dirichlet
minimizer in $U$ of the Mumford-Shah functional.
\end{thm}

This theorem generalizes the result of Example~4 in \cite{Alb-Bou-DM},
where the functions $u_i$ were three distinct constants. The proof is
obtained by the calibration method adapted in \cite{Alb-Bou-DM} to the
functional (\ref{ms}). We construct an explicit calibration for $u$ in a
cylinder $U{\times}\real$, where $U$ is a suitable neighbourhood of $(0,0)$.
The symmetry due to the $2\pi/3$-angles is exploited in the whole construction
of the calibration; in particular, it allows to deduce from the other Euler
conditions that each $u_i$ must be either symmetric or antisymmetric with
respect to the bisecting line of $A_i$ and then, it can be harmonically 
extended to a neighbourhood of the origin, cut by a half-line in the
antisymmetric case. Around the graph of $u$, the calibration is
obtained using the gradient field of a family of harmonic functions, whose
graphs fibrate a neighbourhood of the graph of $u$; this technique reminds the
classical method of the Weierstrass fields, where the minimality of a
candidate $u$ is proved by constructing a suitable slope field, starting from
a family of solutions of the Euler equation, whose graphs fibrate a
neighbourhood of the graph of $u$.

The assumption of $C^2$-regularity for $u_i$ does not seem too restrictive:
indeed, by the regularity results for elliptic problems in non-smooth domains
(see \cite{Gri}), it follows that $u_i$ belongs at least to $C^1(\ol{A_i})$,
since $u_i$ solves the Laplace equation with Neumann boundary
conditions on a sector of angle $2\pi/3$. Moreover,
since $u_i$ is either symmetric or antisymmetric with
respect to the bisecting line of $A_i$, one can see $u_i$ as a solution of 
the Laplace equation on a $\pi/3$-sector with Neumann boundary conditions or
respectively mixed boundary conditions. By the regularity results in
\cite{Gri}, it turns out that in the first case $u_i$ belongs to
$C^2(\ol{A_i})$, while in the second one $u_i$ can be written as 
$u_i(r,\theta)=\td{u}_i(r,\theta)+cr^{3/2}\cos\frac{3}{2}\theta$, with
$\td{u}_i\in C^2(\ol{A_i})$ and $c\in\real$. So, only the function
$r^{3/2}\cos\frac{3}{2}\theta$ is not recovered by our theorem.

The case where $S_u$ is given by three regular curves (not necessarily
rectilinear) meeting at a point with $2\pi/3$-angles, is at the moment an open
problem and it does not seem to be achievable with a plain arrangement of the calibration 
used for the rectilinear case, essentially because of the lack of symmetry
properties.

The paper is organized as follows. In Section~2 we recall the main result of
\cite{Alb-Bou-DM}, while Sections~3 -- 7 are devoted to the proof of
Theorem~\ref{ms}: in Section~3 we construct a calibration $\vp$ in the case 
$u_i$ symmetric and we prove that $\vp$ satisfies conditions (a), (b), (c), and
(e) (see the definition of calibration in Section~2); in Sections~4 and 5 we
show some estimates, which will be useful in Section~6
to prove condition (d); finally, in Section~7 we adapt the calibration to the
antisymmetric case.

\vspace{.5cm}
\noindent {\bf Acknowledgements.}
The author wish to thank Gianni Dal Maso for many interesting discussions on
the subject of this work and for some helpful suggestions about the writing of
this paper.

\section{Preliminary results}

Let $\Omega$ be an open subset of $\real^2$ with a Lipschitz boundary. 
If $u$ is a function in $SBV(\Omega)$, for every $(x_0,y_0)\in \Omega$
one can define
$$u^+(x_0,y_0):=  \aplimsup_{(x,y)\to (x_0,y_0)} u(x,y),
\qquad 
u^-(x_0,y_0):= \apliminf_{(x,y)\to (x_0,y_0)} u(x,y).$$
We recall that $S_u=\{(x,y)\in\Omega:u^-(x,y)<u^+(x,y)\}$ and that
for ${\cal H}^{1}$-a.e. $(x_0,y_0)\in S_u$ there exists a (uniquely defined)
unit vector $\nu_u(x_0,y_0)$ (which is normal to $S_u$ in an approximate
sense) such that
$$\lim_{r\to 0^+} \frac{1}{{\cal
L}^2(B^{\pm}_r( x_0,y_0))}\int_{B^{\pm}_r(x_0,y_0)}
|u(x,y)- u^{\pm}(x_0,y_0)|\,dx\,dy=0, $$
where $B^{\pm}_r(x_0,y_0)$ is the intersection of the ball of radius $r$
centred at $(x_0,y_0)$ 
with the half-plane $\{ (x,y)\in \real^2: 
\pm(x-x_0, y-y_0)\cdot \nu_u(x_0,y_0) \geq 0 \}$.
For more details see \cite{Amb-Fus-Pal}.

For every vector field $\vp:\Omega{\times}\real\to\real^2{\times}\real$ we
define the maps $\vp^{xy}:\Omega{\times}\real\to\real^2$ and
$\vp^z:\Omega{\times}\real\to\real$  by
$$\vp(x,y,z)=(\vp^{xy}(x,y,z),\vp^z(x,y,z)).$$
We shall consider the collection ${\cal F}$ of all piecewise $C^1$ vector
fields $\vp:\Omega{\times}\real\to\real^2{\times}\real$ with the following
property: there exist a finite family $(U_i)_{i\in I}$ of pairwise disjoint
open subsets of $\Omega{\times}\real$ with Lipschitz boundary whose
closures cover $\Omega{\times}\real$, and a family $(\vp_i)_{i\in I}$
of vector fields in $C^1(\ol{U_i}, \real^2{\times}\real)$ such that $\vp$
agrees at any point with one of the $\vp_i$.

Let $u\in SBV(\Omega)$. A {\it calibration\/} for $u$ is a bounded vector
field $\vp\in{\cal F}$ satisfying the following properties:
\begin{description}
	\item[(a)] $\dive\vp=0$ in the sense of distributions in
$\Omega{\times}\real$;
	\item[(b)] $|\vp^{xy}(x,y,z)|^2\leq 4\vp^z(x,y,z)$ at every continuity
point $(x,y,z)$ of $\vp$;
	\item[(c)] $\vp^{xy}(x,y,u(x,y))=2\nabla u(x,y)$ and
$\vp^z(x,y,u(x,y))=|\nabla u(x,y)|^2$ for a.e. $(x,y)\in\Omega\setminus S_u$;
	\item[(d)] $\left| \displaystyle \int_{t_1}^{t_2}\vp^{xy}(x,y,z)\,
dz\right|^2\leq 1$ for every $(x,y)\in\Omega$ and for every $t_1,t_2\in\real$;
	\item[(e)] $\displaystyle \int_{u^-(x,y)}^{u^+(x,y)}\vp^{xy}(x,y,z)\, dz=\nu_u(x,y)$
for ${\cal H}^1$-a.e. $(x,y)\in S_u$.
\end{description}

The following theorem is proved in \cite{Alb-Bou-DM} and \cite{Alb-Bou-DM2}.

\begin{thm}
If there exists a calibration $\vp$ for $u$, then $u$ is a Dirichlet minimizer
of the Mumford-Shah functional (\ref{ms}) in $\Omega$.
\end{thm}

Finally, we present a lemma (proved in \cite{Mor-Mor}), which allows
to construct a divergence free vector field starting from a family of harmonic
functions. 

\begin{lm}\label{folio}
Let $U$ be an open subset of $\real^2$ and $I$, $J$ be two real intervals. 
Let $u: U {\times} J \to I$ be a function of class $C^1$ such that
\begin{itemize}
	\item	$u(\cdot, \cdot \, ; s)$ is harmonic for every $s\in J$;
	\item	there exists a $C^1$ function $t :U { \times } I \to J$
	such that $u(x,y \, ; t(x,y \, ;z ))=z$. 
\end{itemize}
If we define in $U { \times } I$ the vector field 
$$\phi (x,y,z):= (2 \nabla u(x,y \, ; t(x,y \, ;z )), |\nabla u(x,y\, ;
t(x,y\, ;z ))|^2),$$ 
where $\nabla u(x,y\, ; t(x,y\, ;z ))$ denotes the gradient of
$u$ with respect to the variables $x,y$ computed at the point $(x,y\, ; t(x,y\,
;z ))$, then $\phi$ is divergence free in $U{ \times } I$.
\end{lm}

\section{Construction of the calibration}

Let $\{e^x,e^y\}$ be the canonical basis in $\real^2$ and for $i=1,2$ consider
the vectors $\tau_i=(-1/2,(-1)^i\sqrt{3}/2)$, $\nu_i=((-1)^i\sqrt{3}/2,
1/2)$, which are tangent and normal to the set $S_{i-1,i}$ (see Fig.~1). 
As $u_0(0,0)<u_1(0,0)<u_2(0,0)$, there exists an open neighbourhood $U$ of
$(0,0)$ such that the function $u$ belongs to $SBV(U)$, the discontinuity set $S_u$ of
$u$ on $U$ coincides with $\bigcup_{i<j} (S_{i,j}\cap U)$, and the
oriented normal vector $\nu_u$ to $S_u$ is given by 
$$\nu_u(x,y)=\begin{cases}
\none & \text{for $(x,y)\in S_{0,1}$,} \\
\ntwo & \text{for $(x,y)\in S_{1,2}$,} \\
e^y & \text{for $(x,y)\in S_{0,2}$;} \end{cases}$$
by the assumptions on $u_i$, the function $u$ satisfies the Euler conditions
for (\ref{ms}) in $U$. We will construct a local calibration
$\vp=(\vpxi,\vpz):U{\times}\real\to\real^2{\times}\real$ 
for $u$.

\begin{figure}[h]
\begin{center} 
  \includegraphics[height=0.3\textheight,]{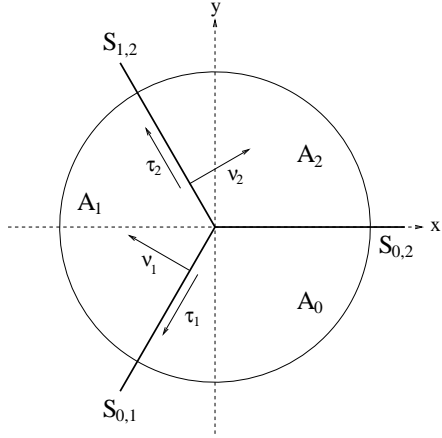}
\caption{the triple junction.}
\end{center}
\end{figure}

Applying Schwarz reflection principle with respect to $S_{0,1}$ and $S_{0,2}$,
the function $u_0$ can be harmonically extended to $U\setminus S_{1,2}$, and
analogously  $u_1$ and $u_2$ can be extended to $U\setminus S_{0,2}$ and
$U\setminus S_{0,1}$, respectively.
By the hypothesis on $u_i$ and by Cauchy-Kowalevski Theorem (see \cite{Joh})
the extension of $u_0$ coincides, up to the sign and to additive constants,
with $u_1$ on $A_1$ and with $u_2$ on $A_2$; analogously, the extension of
$u_1$ coincides, up to the sign and to an additive constant, with $u_2$ on
$A_2$. Since the composition of the three reflections with respect to
$S_{0,1}$, $S_{1,2}$, and $S_{0,2}$ coincides with the reflection with respect
to the bisecting line of the sector $A_0$, by the previous remarks we can
deduce that $u_0$ is either symmetric or antisymmetric with respect to the
bisecting line of $A_0$.

\

We consider first the case $u_0$ symmetric (the antisymmetric case will be
studied in Section~7). Then also $u_1,u_2$ are symmetric with respect to the
bisecting line of $A_1,A_2$, respectively, and the extensions of $u_0,u_1,u_2$
by reflection are well defined and harmonic in the whole set $U$.

In order to define the calibration for $u$, let $\varepsilon >0$,
$l_i\in (u_{i-1}(0,0), u_i(0,0))$ for $i=1,2$, and $\lambda >0$ be
suitable parameters that will be chosen later, and
consider the following subsets of $U { \times } \real$:
$$\begin{array}{rcll}
G_i & := & \{ (x,y,z)\in U{ \times }\real: u_i(x,y)-\eps < z < u_i(x,y)+\eps
\} & \text{for $i=0,1,2$,} \\ 
K_i & := & \{ (x,y,z)\in U{ \times }\real: l_i+\alpha_i(x,y) < z < 
l_i+ 2\lambda+\beta_i(x,y) \} & \text{for $i=1,2$,}\\
H_i & := & \{ (x,y,z)\in U{ \times }\real: l_i +\lambda/2 < z < l_i +3\lambda/2
\} & \text{for $i=1,2$,}
\end{array}$$ 
where $\alpha_i$ and $\beta_i$ are suitable Lipschitz functions such that
$\alpha_i (0,0)= \beta_i(0,0)=0$, which will be defined later. If $\eps$ and 
$\lambda$ are sufficiently small, then for every $i,j$ the sets $G_i$, $K_j$
are nonempty and disjoint, while for every $i$ the set $H_i$ is compactly
contained in $K_i$, provided $U$ is small enough (see Fig.~2).  

\begin{figure}[h!!]
\begin{center} 
  \includegraphics[height=0.7\textheight,]{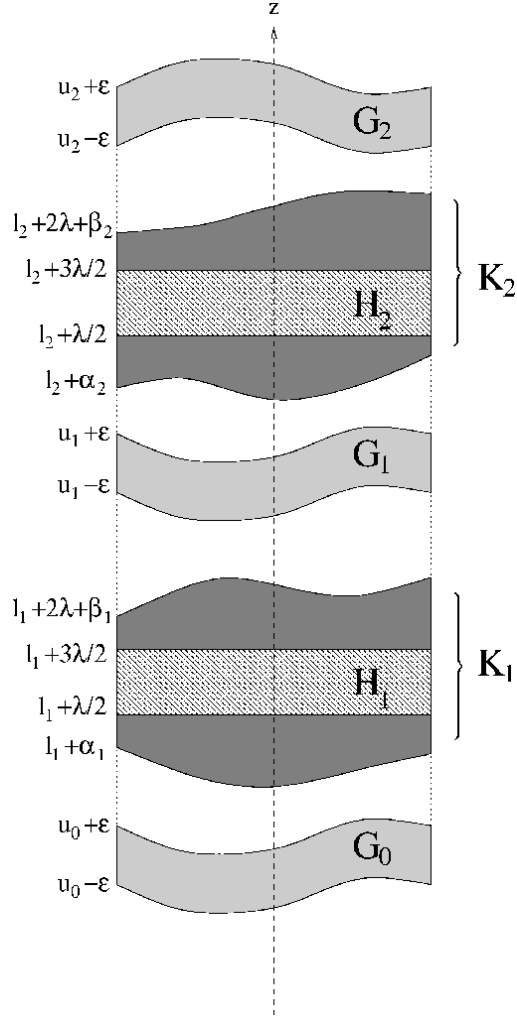}
\caption{section of the sets $G_i, K_i, H_i$ at $x=\text{constant}$.}
\end{center}
\end{figure}

The aim of the definition of the calibration $\vp$ in $G_i$ is
to provide a divergence free vector field satisfying condition (c) and such
that  
$$\begin{array}{c}
\vpxi(s\tau_i,z)\cdot \nu_i > 0 \quad \text{for $u_{i-1}<z<u_{i-1}+\eps$ and
for $u_i-\eps<z<u_i$,} \\
\vpxi(s\tau_i,z)\cdot \nu_i < 0 \quad \text{for $u_{i-1}-\eps<z<u_{i-1}$
and for $u_i<z<u_i+\eps$,}
\end{array}$$
for $i=1,2$ and $s\geq 0$, and analogously
$$\begin{array}{c}
\vpxi(s,0,z)\cdot e^y > 0 \quad \text{for $u_0<z<u_0+\eps$
and for $u_2-\eps<z<u_2$,} \\
\vpxi(s,0,z)\cdot e^y < 0 \quad \text{for $u_0-\eps<z<u_0$ and
for $u_2<z<u_2+\eps$;}
\end{array}$$ 
these properties are crucial in order to obtain (d) and (e) simultaneously.
Such a field can be obtained by applying the technique shown in
Lemma~\ref{folio}, starting from the family of harmonic functions $u_i+tv_i$,
where we choose as $v_i$ the linear functions defined by
$$v_0(x,y):= \ttwo\!\cdot\!(x,y) + \eps, \
v_1(x,y):= \tzer\!\cdot\!(x,y) +\eps, \
v_2(x,y):= \tone\!\cdot\!(x,y) +\eps.$$
So for every $(x,y,z)\in G_i$, $i=0,1,2$, we define the vector $\vp(x,y,z)$ as
$$\displaystyle \left( 2 \nabla u_i + 2\,
\frac{z-u_i(x,y)}{v_i(x,y)}\nabla v_i, \left| \nabla u_i
+ \frac{z-u_i(x,y)}{v_i(x,y)}\nabla v_i \right|^2 \right).$$

The r\^ole of $K_i$ is to give the exact contribution to the integral in (e).
In order to annihilate the tangential contribution on $S_u$ due to the
choice of the field in $G_i$, we insert in $K_i$ the region $H_i$ and 
for every $(x,y,z)\in H_i$, $i=1,2$, we define $\vp(x,y,z)$ as
$$\displaystyle \left(-\frac{2\eps}{\lambda}
\left(\nabla u_{i-1}+\nabla u_i\right) , \mu \right)$$
where $\mu$ is a positive constant which will be suitably chosen later.
By the harmonicity of $u_i$ this field is divergence free and, as
$\df_{\nu}u_i=0$ on $S_u$ for every $i$, its horizontal component is purely
tangential on $S_u$. So, it remains to correct only the normal contribution to
the integral in (e)  due to the field in $G_i$. To realize this purpose on the
two segments $S_{i-1,i}$, $i=1,2$, we could require that
$\alpha_i(s\tau_i)=\beta_i(s\tau_i)=0$ for every $s\geq 0$ (see the
definition of $K_i$) and define $\vp(x,y,z)$ for $(x,y,z)\in K_i\setminus
\ol{H_i}$ as 
\begin{equation}\label{expl} \displaystyle \left(
\frac{1}{\lambda} g\left(\tau_i\!\cdot\!(x,y)\right)\nu_i, \mu \right),
\end{equation} 
where $g$ is a function of real variable chosen in such a way
that (e) is satisfied for $(x,y)\in S_{i-1,i}$, i.e.,
$$g(t):= 1 - \sqrt{3}\frac{\eps^2}{v_0(t,0)} \qquad
\forall t\in \real,$$
as we will see later in (\ref{ei}). Note that the two-dimensional field
$g\left(\tau_i\!\cdot\!(x,y)\right)\nu_i$ is divergence free, since it is with
respect to the orthonormal basis $\{\tau_i,\nu_i\}$, hence $\vp$ is divergence
free in $K_i\setminus\ol{H_i}$; moreover, since $\vpz\equiv \mu$ on $K_i$, the
normal component of $\vp$ is continuous across the boundary of $H_i$, so that
$\vp$ turns out to be divergence free in the sense  of distributions
in the whole set $K_i$. Actually it is crucial to add a component along the
direction $\tau_i$ to the field in (\ref{expl}) in order to make (d) true, as
it will be clear in the proof of Step~2 (see Section~5); this component has to
be chosen in such a way that it is zero on $S_{i-1,i}$ (so that (e) remains
valid on these segments) and that it depends only on $\nu_i\!\cdot\!(x,y)$ (so
that the field remains divergence free). Therefore we replace in (\ref{expl})
the vector $g\left(\tau_i\!\cdot \!(x,y)\right)\nu_i$ by 
\begin{equation}\label{phi}
\phi_i(x,y):= (-1)^{i+1}f\left(\nu_i\!\cdot\!(x,y)\right)\tau_i +
g\left(\tau_i\!\cdot\!(x,y)\right)\nu_i, 
\end{equation}
where $f$ is an even smooth function of real variable such that $f(0)=0$
and which will be chosen later in a suitable way (see (\ref{condf})).
From this definition it follows that 
\begin{equation}\label{phii} \phi_2^{x}(x,y)= -\phi_1^{x}(x, -y), \qquad
\phi_2^{y}(x,y)= \phi_1^{y}(x, -y),
\end{equation}
so that
$$\phi_1(x,0)+\phi_2(x,0)=2\phi_1^{y}(x,0)\nzer,$$
i.e., if we assume that $\alpha_i(x,0)=\beta_i(x,0)$ for every $x\geq 0$, the
contribution given by the fields (\ref{phi}) to the integral in (e) computed
at a point of $S_{0,2}$ is purely normal, as required in (e), but its modulus
is in general different from what we need to obtain exactly the normal
vector $\nzer$. In order to correct it, we multiply $\phi_i$ by a function
$\sigma_i$ which is first defined on $S_{i-1,i}\cup S_{0,2}$ (more precisely,
$\sigma_i$ is taken equal to $1$ on $S_{i-1,i}$ and to the correcting factor
on $S_{0,2}$); then, we extend it to a neighbourhood of $(0,0)$ by assuming
$\sigma_i$ constant along the  integral curves of $\phi_i$, so that
$\sigma_i\phi_i$ remains divergence free. 

The integral curves of $\phi_i$ can be represented as the curves
$\{(x,y)\in U:y=\psi_i(x,s)\}$, where $\psi_i(x,s)$ is the solution of the
problem \begin{equation}\label{psidef}
\begin{cases}
\df_{x}\psi_i(x, s)\phi_i^{x}(x, \psi_i(x,s)) -
\phi_i^{y}(x, \psi_i(x,s))=0, \\
\psi_i(s,s)=0,
\end{cases}
\end{equation}
which is defined in a sufficiently small neighbourhood of $(0,0)$. By applying
the Implicit Function Theorem, it is easy to see that if $U$ is small enough,
then there exists a unique smooth function $h_i$ defined in $U$ such that 
\begin{equation}\label{hdef}
h_i(0,0)=0, \qquad \psi_i(x,h_i(x,y))=y.
\end{equation}
Note that the curve $\{(x,y)\in U:h_i(x,y)=s\}$ coincides with the integral curve
$\{(x,y)\in U:y=\psi_i(x,s)\}$ and that $(h_i(x,y),0)$ gives the intersection
point of the integral curve passing through $(x,y)$ with the $x$-axis; in
other words, the level lines of $h_i$ provide a different representation of the
integral curves of $\phi_i$ in terms of their intersection point with the
$x$-axis.
\begin{figure}[h]
\begin{center} 
	\includegraphics[height=0.3\textheight,]{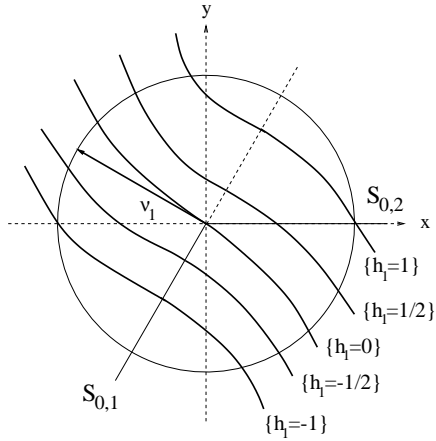}
\caption{integral curves of the field $\phi_1$.}
\end{center}
\end{figure}

We state here some properties of $h_i$ and $\psi_i$ for further
references. Since $\psi_i(s,s)=0$, we have that 
\begin{equation}\label{sus}
h_i(s,0)=s
\end{equation}
for every $s$ such that $(s,0)\in U$. By (\ref{psidef}) and by differentiating
the initial condition in (\ref{psidef}) with respect to $s$, we obtain
\begin{equation}\label{psider}
\df_x\psi_i(0,0)=\frac{\phi_i^y(0,0)}{\phi_i^x(0,0)}=\frac{\nu_i^y}{\nu_i^x}=
\frac{(-1)^i}{\sqrt{3}}, 
\qquad \df_s\psi_i(0,0)=-\df_x\psi_i(0,0)=\frac{(-1)^{i+1}}{\sqrt{3}}.
\end{equation}
By differentiating the equation in (\ref{psidef}) with respect to $x$ and to
$s$, and by using (\ref{phi}), it is easy to see that
\begin{equation}\label{psider2}
\df^2_x\psi_i(0,0)=\df^2_{xs}\psi_i(0,0)=0,
\end{equation}
while by differentiating twice with respect to $s$ the initial condition
$\psi_i(s,s)=0$, we obtain that 
\begin{equation}\label{psider2bis}
\df^2_s\psi_i(0,0)=-2\df^2_{xs}\psi_i(0,0)=0.
\end{equation}
By (\ref{psider}) and (\ref{psider2}), the curve $\{h_i=0\}$ (which coincides
with $\{y=\psi_i(x,0)\}$) is tangent to $\nu_i$ at $0$, which may be an
inflection point. Moreover, since $\df_x\psi_i(0,0)\neq 0$, by continuity
the function $\psi_i(\cdot,s)$ is strictly monotone in a small neighbourhood of
$0$ for $s$ sufficiently small; by this fact and by comparing the values of
the function $\psi_i(\cdot, h_i(s\tau_i))$ at the points $h_i(s\tau_i)$ and
$s\tau_i^x$, it is easy to see that
\begin{equation}\label{susi} 
h_i(s\tau_i)\leq 0 
\end{equation}
for every $s\geq 0$ such that $s\tau_i\in U$, provided $U$ is small enough.
Remark that by (\ref{sus}) and (\ref{susi}) it follows that the segment
$S_{0,2}$ is all contained in the region $\{h_i\geq 0\}$, while
$S_{i-1,i}$ in the region $\{h_i\leq 0\}$.

At last, we set
$$\sigma_i(x,y):= 
\begin{cases}
1 & \text{if $h_i(x,y)\leq 0$,} \\
\displaystyle \frac{g(h_i(x,y))}{2\phi_i^{y}(h_i(x,y),0)} &
\text{if $h_i(x,y)>0$;}
\end{cases}$$
since by definition $\phi_i^y(0,0)=g(0)\nu_i^y=g(0)/2$, the function 
$\sigma_i$ is continuous across the curve $\{h_i=0\}$. Moreover, 
remark that from (\ref{phii}) it
follows that $\psi_2(x, s)=-\psi_1(x,s)$,
$h_2(x,y)=h_1(x,-y)$, and then
\begin{equation}\label{sigmai}
\sigma_2(x,y)=\sigma_1(x,-y).
\end{equation}
For every $(x,y,z)\in K_i\setminus\ol{H_i}$, $i=1,2$, we define $\vp(x,y,z)$ as
$$\displaystyle \left( \frac{1}{\lambda}\sigma_i(x,y) \phi_i(x,y), \mu
\right).$$

In the remaining regions of transition it is convenient to take $\vp$ purely
vertical. In order to make $\vp$ divergence free in the whole set
$U{\times}\real$, we need the normal component of $\varphi$ to be 
continuous across the boundary of $G_i$ and $H_i$. To guarantee this
continuity across $\df G_i$, we are forced to take as third component
of $\vp$ the function  \begin{equation}
\omega(x,y,z):=
\begin{cases}
\displaystyle \frac{\eps^2}{v_0^2(x,y)}-|\nabla u_0|^2 & \text{for
$z<l_1+\lambda$,} \\
\\
\displaystyle \frac{\eps^2}{v_1^2(x,y)}-|\nabla u_1|^2 & \text{for
$l_1+\lambda
\leq z<l_2+\lambda$,} \\
\\
\displaystyle \frac{\eps^2}{v_2^2(x,y)}-|\nabla u_2|^2 & \text{for
$z\geq l_2+\lambda$.}
\end{cases}
\end{equation}

Finally, we define the functions $\alpha_i,\beta_i$ in such a way that the
normal component of $\vp$ turns out to be continuous also across the
boundary of $K_i$; more precisely, for $i=1,2$ we choose $\alpha_i$ as
the solution of the Cauchy problem
$$\begin{cases}
\displaystyle
\frac{1}{\lambda}\sigma_i(x,y)\phi_i(x,y)
\cdot\nabla \alpha_i(x,y) -\mu =
-\frac{\eps^2}{v_{i-1}^2(x,y)} +|\nabla u_{i-1}(x,y)|^2, \\
\alpha_i(s\tau_i)=0, \ \alpha_i(s,0)=0 \quad \text{for $s\geq 0$,} 
\end{cases}$$
while $\beta_i$ as the solution of
$$\begin{cases}
\displaystyle
\frac{1}{\lambda}\sigma_i(x,y)\phi_i(x,y)
\cdot\nabla \beta_i(x,y) -\mu =
-\frac{\eps^2}{v_i^2(x,y)} +|\nabla u_i(x,y)|^2, \\
\beta_i(s\tau_i)=0,\ \beta_i(s,0)=0 \quad \text{for $s\geq 0$.} 
\end{cases}$$
Since $\sigma_i$ is not $C^1$ near the curve $\{h_i=0\}$, we cannot expect a
$C^1$ solution. Nevertheless, if $U$ is small enough, then $\alpha_i,\beta_i$
are Lipschitz function defined in $U$, and the possible discontinuity points of
$\nabla\alpha_i, \nabla\beta_i$ concentrate only on the curve $\{h_i=0\}$;
indeed, if $U$ is sufficiently small, the Cauchy problems
\begin{equation}\label{aux1}
\begin{cases}
\displaystyle\frac{1}{\lambda}\phi_i(x,y)
\cdot\nabla \tilde\alpha_i(x,y) -\mu =
-\frac{\eps^2}{v_{i-1}^2(x,y)} +|\nabla u_{i-1}(x,y)|^2, \\
\tilde\alpha_i(s\tau_i)=0 \quad (s\in\real),
\end{cases}
\end{equation}
and
\begin{equation}\label{aux2}
\begin{cases}
\displaystyle\frac{g(h_i(x,y))}{2\lambda\phi_i^{y}(h_i(x,y),0)}\phi_i(x,y)
\cdot\nabla \hat\alpha_i(x,y) -\mu =
-\frac{\eps^2}{v_{i-1}^2(x,y)} +|\nabla u_{i-1}(x,y)|^2, \\
\hat\alpha_i(s,0)=0 \quad (s\in\real),
\end{cases}
\end{equation}
admit a unique solution $\tilde\alpha_i,\hat\alpha_i\in C^{\infty}(U)$, since the lines
$\{s\tau_i:s\in\real\}$ and $\{(s,0):s\in\real\}$ are not characteristic for
these equations. Since the curve $\{h_i=0\}$, which coincides with the curve
$\{y=\psi_i(x,0)\}$, is a characteristic line of both equations (\ref{aux1})
and (\ref{aux2}) (use (\ref{psidef}) and $g(0)/(2\lambda\phi_i^y(0,0))=1$),
the functions $\tilde\alpha_i,\hat\alpha_i$ assume
the same value on the curve $\{h_i=0\}$. So, $\alpha_i$ can be regarded as the
function defined by $$\alpha_i(x,y):= 
\begin{cases}
\tilde\alpha_i(x,y) & \text{if $h_i(x,y)\leq 0$,} \\
\hat\alpha_i(x,y) & \text{if $h_i(x,y)>0$,}
\end{cases}$$
and therefore $\alpha_i$ is $C^{\infty}$ in $U\setminus\{h_i=0\}$, and all
derivatives of $\alpha_i$ have finite limits on both sides of $\{h_i=0\}$.
The same argument works for $\beta_i$.

The complete definition of the field is therefore the following: for every 
$(x,y,z)\in U {\times } \real$, the vector 
$\vp(x,y,z)=(\vpxi, \vpz)(x,y,z)\in \real^2 { \times } \real$
is given by
%
%
$$
\begin{cases}
\displaystyle \left( 2 \nabla u_i + 2\, \frac{z-u_i(x,y)}{v_i(x,y)}\nabla
v_i, \left| \nabla u_i
+ \frac{z-u_i(x,y)}{v_i(x,y)}\nabla v_i \right|^2 \right) 
& \text{in $G_i$ \ $(i=0,1,2)$,} \\
\\
\displaystyle \left( \frac{1}{\lambda}\sigma_i(x,y) \phi_i(x,y), \mu \right) &
\text{in  $K_i \setminus \overline{H_i}$ \ $(i=1,2)$,} \\
\\
\displaystyle \left(-\frac{2\eps}{\lambda}
\left(\nabla u_{i-1}+\nabla u_i\right) , \mu \right) & \text{in $H_i$ \
$(i=1,2)$,}\\ \\
\displaystyle \left( 0, \omega(x,y,z)\right) & \text{otherwise.} 
\end{cases}$$

By construction conditions (a) and (c) are satisfied.

Condition (b) is trivial in $G_i$ for all $i$.

Since $\nabla u_i(0,0)=0$ for all $i$ (this fact easily follows by the assumptions
on the regularity of $u_i$ and by the Euler conditions), we have that 
$$\frac{\eps^2}{v^2_i(0,0)} - |\nabla u_i(0,0)|^2 = 1 >0;$$
then, if $U$ is small enough,  
$$\frac{\eps^2}{v^2_i(x,y)} - |\nabla u_i(x,y)|^2 >0$$
for every $(x,y)\in U$ and for every $i=0,1,2$, and so $\omega$ is always
positive.

Arguing in a similar way, if we impose that
$\mu > 1/(4\lambda^2)$, 
condition (b) holds in $K_i$, 
provided $U$ is sufficiently small.

By direct computations we find that for every $(x,y)\in U$
\begin{align}\label{01}
\int_{u_{i-1}}^{u_i}\vpxi\ dz & =  
\frac{\eps^2}{v_{i-1}}\nabla v_{i-1} -
\frac{\eps^2}{v_i}\nabla v_i
+\frac{1}{\lambda}(\beta_i - \alpha_i+ \lambda)\sigma_i
\phi_i, \\ \intertext{for $i=1,2$, while}
\int_{u_0}^{u_2}\vpxi \label{02}
\ dz & = 
\frac{\eps^2}{v_0}\nabla v_0 -
\frac{\eps^2}{v_2}\nabla v_2
+\frac{1}{\lambda}\sum_{i=1}^2(\beta_i - \alpha_i + \lambda)\sigma_i
\phi_i.
\end{align}
Note that for $i=1,2$
\begin{eqnarray}\label{veq}
& v_{i-1}(s\tau_i)=v_i(s\tau_i)=v_0(s,0)=\displaystyle -\frac{s}{2}+\eps & 
\forall s\in\real,  \\ \vs & \nabla v_{i-1}(x,y)-\nabla v_i(x,y)=\sqrt{3}\nu_i
\label{nbl-} &  \forall (x,y)\in U.
\end{eqnarray}
As $h_i(s\tau_i)\leq 0$ for every $s\geq 0$ by (\ref{susi}), we have that
$\sigma_i(s\tau_i)=1$ for every $s\geq 0$, while by definition
$\alpha_i(s\tau_i)=\beta_i(s\tau_i)=0$. From these facts, (\ref{01}),
(\ref{veq}), (\ref{nbl-}), and the definition of $\phi_i$, we obtain 
\begin{equation}\label{ei}
\int_{u_{i-1}(s\tau_i)}^{u_i(s\tau_i)}\vpxi
(s\tau_i,z)\ dz  = \sqrt{3}\frac{\eps^2}{v_0(s,0)}
\nu_i +(-1)^{i+1}f(0)\tau_i + g(s)\nu_i  =\nu_i,
\end{equation}
where the last equality follows from the definition of $g$ and the fact that
$f(0)=0$. Analogously, by the equalities
\begin{eqnarray}\label{veq0}
& v_0(s,0)=v_2(s,0) & \forall s\in \real, \\
& \nabla v_0(x,y)-\nabla v_2(x,y)=\sqrt{3}\nzer & \forall (x,y)\in
U,\label{nbl-0} \end{eqnarray}
by the definition of $\alpha_i$ and $\beta_i$, and by (\ref{phii}),
(\ref{sigmai}), (\ref{02}), we have 
\begin{eqnarray}
\int_{u_0(s,0)}^{u_2(s,0)}\vpxi
(s,0,z)\ dz & = & \sqrt{3}\frac{\eps^2}{v_0(s,0)}\nzer 
+ 2\sigma_1(s,0)\phi_1^y(s,0)\nzer \nonumber \\
& = & \sqrt{3}\frac{\eps^2}{v_0(s,0)}\nzer
+ g(s)\nzer =\nzer, \label{e0}
\end{eqnarray}
where the two last equalities follow from (\ref{sus}) and from the definition
of $\sigma_1$ and  $g$. So condition (e) is satisfied.

The proof of condition (d) will be split in the next three sections: in
Section~4 we prove that condition (d) holds if $t_1$ and $t_2$ belong
to suitable neighbourhoods of $u_{i-1}(0,0)$ and $u_i(0,0)$, respectively;
then, in Section~5 we prove condition (d) for $t_1$ and $t_2$ belonging 
to suitable neighbourhoods of $u_0(0,0)$ and $u_2(0,0)$,
respectively; finally, in Section~6, by a continuity argument we show that
condition (d) is true in all other cases.

\section{Estimates for $t_1$ and $t_2$ near $u_{i-1}$ and $u_i$}

For $(x,y)\in U$ and $t_1,t_2\in\real$, we set
\begin{equation}\label{intgr}
I(x,y,t_1,t_2):= \int_{t_1}^{t_2} \vpxi(x,y,z)\,dz
\end{equation}
and we denote its absolute value by $\rho$. In this section, we will show
that $\rho(x,y,t_1,t_2)\leq 1$ in a neighbourhood of the point
$(0,0,u_{i-1}(0,0),u_i(0,0))$ for $i=1,2$, so that the following step 
will be proved.

\bstep{1}{For a suitable choice of the parameter $\varepsilon$, 
there exists $\delta >0$ such that condition (d)
holds for $|t_1-u_{i-1}(0,0)|<\delta$, $|t_2-u_i(0,0)|<\delta$ with $i=1,2$,
provided $U$ is small enough.}

\noindent
Note that $\rho$ is a continuous
function, but its derivatives with respect to $x,y$ may be discontinuous at
the points $(x,y,t_1,t_2)$ such that $h_1(x,y)=0$ or $h_2(x,y)=0$; indeed, the
curve $\{h_{i}=0\}$ is the boundary of the different regions of definition of
the functions $\sigma_i$, $\alpha_i$, and $\beta_i$, whose derivatives may
present therefore some discontinuities. Nevertheless, if we set
$N_i:=\{(x,y)\in U:h_i(x,y)<0\}$ and $P_i:=\{(x,y)\in U:h_i(x,y)>0\}$, the
restrictions of $\sigma_{i}$, $\alpha_{i}$, and  $\beta_i$ to the sets $N_i$
and $P_i$ can be extended up to the boundary $\{h_i=0\}$ as $C^{\infty}$
functions; so, along the curve $\{h_i=0\}$ the traces of the derivatives of
$\sigma_{i}$, $\alpha_{i}$, and $\beta_i$ are defined. Then, also the traces
of the derivatives of $\rho$ with respect to $x,y$ are defined at the points
$(x,y,t_1,t_2)$ with $h_1(x,y)=0$ or $h_2(x,y)=0$.

\begin{figure}[h]
\begin{center} 
	\includegraphics[height=0.3\textheight,]{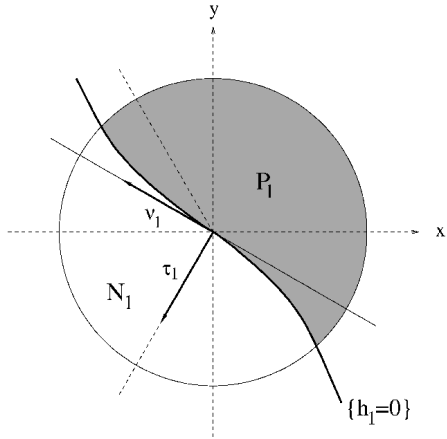}
\caption{the regions $P_1$ and $N_1$.}
\end{center}
\end{figure}

Since we want to study the behaviour of $\rho$ in a neighbourhood of
$(0,0,u_{i-1}(0,0), u_i(0,0))$, we can suppose $|t_1-u_{i-1}(0,0)|\leq\eps$ and
$|t_2-u_i(0,0)|\leq\eps$, so that the possible discontinuities of the
derivatives of $\rho$ concentrate only on the curve $\{ h_i=0\}$. We study
separately the two regions $N_i$ and $P_i$.

\vspace{.1cm}

Consider first the case $(x,y)\in \ol{N_i}$, which is the region containing
$S_{i-1,i}$. We will study the derivatives of $\rho$ at the points of the form 
$$q_i(s):=(s\tau_i, u_{i-1}(s\tau_i),u_i(s\tau_i)), \qquad s\geq 0.$$
We have already shown (condition (e)) that
$\rho(q_i(s)) =1$ for every $s\geq 0$; we want
to prove that 
\begin{equation}\label{null}
\nabla\rho(q_i(s))=0 \qquad \forall s\geq 0
\end{equation}
(where now $\nabla$ denotes the gradient with respect to $x,y,t_1,t_2$) and
that the Hessian matrix of $\rho$ with respect to $\nu_i,t_1,t_2$ is negative
definite at $q_i(0)$. 

Let $I^{\tau_i}$ and $I^{\nu_i}$ be the
components of the integral in (\ref{intgr}) along the directions $\tau_i$ and
$\nu_i$, respectively. Since by definition 
$$\rho(x,y,t_1,t_2)=[(I^{\tau_i}(x,y,t_1,t_2))^2+(I^{\nu_i}(x,y,t_1,t_2))^2]^{1/2},$$
the gradient of $\rho$ is given by
\begin{equation}\label{mulgrad}
\nabla \rho = \frac{1}{\rho}(I^{\tau_i}\nabla I^{\tau_i} + I^{\nu_i} \nabla
I^{\nu_i}).
\end{equation}
Note that (\ref{ei}) implies that 
\begin{equation}\label{cpei}
I^{\tau_i}(q_i(s))=0 \qquad \text{and} \qquad I^{\nu_i}(q_i(s))=1 \qquad
\forall s\geq 0,
\end{equation}
hence
\begin{equation}\label{equal}  
\nabla \rho (q_i(s))=\nabla I^{\nu_i}(q_i(s)).
\end{equation} 
By the definition of $\vp$ in $G_i$ and by (\ref{01}) we can compute
explicitly the expression of $I^{\nu_i}$ at $(x,y,t_1,t_2)$:
\begin{multline}\label{expr} 
I^{\nu_i} = -2(t_1-u_{i-1})\dfi u_{i-1}+2(t_2-u_i)\dfi u_i+
\frac{1}{\lambda}(\beta_i -\alpha_i
+\lambda)\sigma_i\phi_i^{\nu_i} \\
+\frac{\sqrt{3}}{2v_{i-1}}(\eps^2-(t_1-u_{i-1})^2)+\frac{\sqrt{3}}{2v_i}(\eps^2-(t_2-u_i)^2),
\end{multline} 
where
\begin{equation}\label{prime}
\phi_i^{\tau_i}(x,y)=(-1)^{i+1}f(\nu_i\!\cdot\!(x,y)) \qquad \text{and} \qquad
\phi_i^{\nu_i}(x,y)=g(\tau_i\!\cdot\!(x,y)).
\end{equation}
By differentiating (\ref{expr}) with respect to the direction $\nu_i$ we
obtain  \begin{multline}\label{fnd}
\dfi I^{\nu_i} = 2(\dfi u_{i-1})^2-2(\dfi u_i)^2
-2(t_1-u_{i-1})\dfii u_{i-1}+2(t_2-u_i)\dfii u_i \\
+\frac{1}{\lambda}\dfi (\beta_i-
\alpha_i)\sigma_i\phi_i^{\nu_i} 
+ \frac{1}{\lambda}(\beta_i
-\alpha_i +\lambda)(\dfi\sigma_i \phi_i^{\nu_i} +
\sigma_i\dfi\phi_i^{\nu_i}) \\
-\frac{3}{4v_{i-1}^2}(\eps^2-(t_1-u_{i-1})^2)
+\frac{3}{4v_i^2}(\eps^2-(t_2-u_i)^2) \\
+\frac{\sqrt{3}}{v_{i-1}}(t_1-u_{i-1})\dfi
u_{i-1}+\frac{\sqrt{3}}{v_i}(t_2-u_i)\dfi u_i.
\end{multline}
By the Euler conditions, $\dfi u_{i-1}(s\tau_i)=\dfi u_i(s\tau_i)=0$
for every $s\geq 0$. Moreover, since $|\nabla u_{i-1}|=|\nabla u_i|$ on $U$
(see the remark at the beginning of the proof),
in the region $\ol{N_i}$ the function
$\beta_i-\alpha_i$  coincides with the solution $\xi_i$ of the problem 
\begin{equation}\label{dfr}
\begin{cases}
\displaystyle
\frac{1}{\lambda}\phi_i^{\tau_i}\dft\xi_i+\frac{1}{\lambda}\phi_i^{\nu_i}
\dfi\xi_i= \frac{\eps^2}{v_{i-1}^2}-\frac{\eps^2}{v_i^2} , \\
\xi_i(s\tau_i)= 0 \quad (s\geq 0). 
\end{cases}
\end{equation}
As $\dft\xi_i(s\tau_i)=0$ and $v_{i-1}(s\tau_i)=v_i(s\tau_i)$ for every $s\geq
0$ (see (\ref{veq})), we have that
\begin{equation}\label{null1}
\dfi (\beta_i-\alpha_i)(s\tau_i)=\dfi \xi_i (s\tau_i)=0.
\end{equation}
By definition $\dfi \phi_i^{\nu_i}\equiv 0$ and $\sigma_i(x,y)=1$ 
for every $(x,y)\in \ol{N_i}$; using these remarks and the first equality in
(\ref{veq}), we can deduce that  
\begin{equation}\label{null0}
\dfi I^{\nu_i}(q_i(s))=0
\end{equation}
for every $s>0$, and the equality holds also for the trace of $\dfi
I^{\nu_i}$ at $q_i(0)$.
Since the derivatives of $I^{\nu_i}$ with respect to $t_1$ and $t_2$ are given
by \begin{equation}\label{26}
\df_{t_1}I^{\nu_i}=-2\dfi u_{i-1}-\frac{\sqrt{3}}{v_{i-1}}(t_1-u_{i-1}),
\qquad
\df_{t_2}I^{\nu_i}=2\dfi u_i-\frac{\sqrt{3}}{v_i}(t_2-u_1),
\end{equation}
by the Euler conditions it follows that 
\begin{equation}\label{2600}
\df_{t_1}I^{\nu_i}(q_i(s))=
\df_{t_2}I^{\nu_i}(q_i(s))=0.
\end{equation}
As $I^{\nu_i}(q_i(s))=1$ for every $s\geq 0$, equalities (\ref{2600}) imply
that $\df_{\tau_i}I^{\nu_i}(q_i(s))=0$. By this fact,
(\ref{equal}), (\ref{null0}), and (\ref{2600}),
equality (\ref{null}) is proved.

Now we need to compute the trace of the Hessian matrix of $\rho$ with respect 
to $\nu_i,t_1,t_2$ at the point $q_i(0)$; using (\ref{cpei})
(\ref{null0}), (\ref{2600}) and (\ref{null}), the
Hessian matrix at $q_i(0)$ reduces to 
\begin{equation}\label{grads}
\hess \rho(q_i(0)) = [\nablav I^{\tau_i} 
\otimes\nablav I^{\tau_i}+\hess I^{\nu_i}](q_i(0)),  \end{equation}
where $\nablav$ denotes the gradient with respect to $\nu_i,t_1,t_2$ and
$\otimes$ the tensor product. As before, we know the explicit expression of
$I^{\tau_i}$: \begin{multline}\label{exprt}
I^{\tau_i} = -2(t_1-u_{i-1})\dft u_{i-1}+2(t_2-u_i)\dft u_i+
\frac{1}{\lambda}(\beta_i -\alpha_i
+\lambda)\sigma_i\phi_i^{\tau_i} \\
-\frac{1}{2v_{i-1}}(\eps^2-(t_1-u_{i-1})^2)+\frac{1}{2v_i}(\eps^2-(t_2-u_i)^2),
\end{multline}
hence, using the Euler conditions, (\ref{null1}), and the fact that
$\sigma_i\equiv 1$ in $\ol{N_i}$, it results that
\begin{equation}\label{tg}
\dfi I^{\tau_i}(q_i(0))  =  \frac{1}{2}\dfi v_{i-1}(0,0) -\frac{1}{2}
\dfi v_i(0,0) + \dfi \phi_i^{\tau_i}(0,0)
=  \frac{\sqrt{3}}{2}, 
\end{equation}
where the last equality follows by (\ref{nbl-}) and by the equality
\begin{equation}\label{dntn}
\dfi \phi_i^{\tau_i}(0)=(-1)^{i+1}f'(0)=0.
\end{equation}
By differentiating (\ref{fnd}) and by using the Euler conditions,
(\ref{null1}), the constancy of $\sigma_i$ in $\ol{N_i}$, and the fact that
$\dfii \phi_i^{\nu_i}\equiv 0$, we have
\begin{equation}\label{nm} 
\dfii I^{\nu_i}(q_i(0))  = \frac{1}{\lambda}\phi_i^{\nu_i}(0,0)\dfii (\beta_i
-\alpha_i)(0,0)+
\frac{3}{2\eps}\dfi v_{i-1}(0,0)-
\frac{3}{2\eps}\dfi v_i(0,0)=
-\frac{\sqrt{3}}{2\eps},
\end{equation}
where the last equality follows from
\begin{equation}\label{star}
\frac{1}{\lambda}\phi_i^{\nu_i}(0,0)\dfii (\beta_i-\alpha_i)(0,0)=
-\frac{2\sqrt{3}}{\eps},
\end{equation}
which can be obtained by differentiating (\ref{dfr}). 
Using (\ref{grads}), (\ref{tg}), and (\ref{nm}), we obtain that
\begin{equation}\label{final1}
\dfii \rho(q_i(0))=[\dfi I^{\tau_i}(q_i(0))]^2 + \dfii I^{\nu_i}(q_i(0))= \frac{3}{4}-
\frac{\sqrt{3}}{2\eps} <0, 
\end{equation}
provided $\eps$ is sufficiently small. 
Since $\df_{t_1}I^{\tau_i}(q_i(0))=0$ (this can be
easily proved using the fact that $\nabla u_{i-1}(0,0)=\nabla u_i(0,0)=0$), by
(\ref{grads}) we have that
$$\dsns \rho(q_i(0))=\dsns I^{\nu_i}(q_i(0)) , \qquad
\df^2_{t_1} \rho(q_i(0))=\df^2_{t_1} I^{\nu_i} (q_i(0)).$$
By differentiating (\ref{26}) and by using the Euler conditions, it turns out
that $$\dsns I^{\nu_i}(q_i(0))=-2\dfii u_{i-1}(0,0), \qquad \df^2_{t_1}
I^{\nu_i} (q_i(0))=-\frac{\sqrt{3}}{\eps},$$ 
so that
$$\det\left( \begin{array}{cc}
\dfii \rho & \dsns \rho \vspace{.1cm}\\
\dsns \rho & \df^2_{t_1} \rho
\end{array}\right)(q_i(0))= \frac{3}{2\eps^2}\left(
1-\frac{\sqrt{3}}{2}\eps\right) -4(\dfii
u_{i-1}(0,0))^2.$$
Arguing in a similar way, one can find that
$$\dsnt \rho(q_i(0))= 2\dfii u_i(0,0), \quad 
\df^2_{t_2} \rho(q_i(0)) = -\frac{\sqrt{3}}{\eps}, \quad
\dsst \rho(q_i(0)) =0,$$
so that
$$\det\hess \rho(q_i(0))= -\frac{3\sqrt{3}}{2\eps^3}
\left(1 -\frac{\sqrt{3}}{2}\eps\right)
+\frac{4\sqrt{3}}{\eps}[(\dfii u_{i-1}(0,0))^2+(\dfii u_i(0,0))^2].$$
Since for $\eps$ sufficiently small it results that
\begin{equation}\label{final100}
\det\left( \begin{array}{cc} \dfii \rho & \dsns \rho \vspace{.1cm}\\
\dsns \rho & \df^{2}_{t_1}\rho \end{array}\right) (q_{i}(0))>0, \qquad
\det\hess\rho(q_{i}(0))<0,
\end{equation}
then, by (\ref{final1}) and (\ref{final100}) the Hessian matrix of 
$\rho$ at $q_{i}(0)$ is negative definite.

At this point we have all the ingredients we need in order to compare the
value of $\rho$ on $S_{i-1,i}$ with its value at a point $(x,y,t_1,t_2)$ for
$(x,y)\in\ol{N_i}$ and $|t_1-u_{i-1}(0,0)|\leq\eps$, $|t_2-u_{i}(0,0)|\leq\eps$.

Remark that since the curve $\{ h_i=0 \}$ may have an inflection point
at the origin, the set $\ol{N_i}$ might be not convex. 
If the segment joining $(x,y)$ with its orthogonal projection on $S_{i-1,i}$ 
(which is a point of the form $s\tau_i$ with $s\geq 0$) is all contained in
$\ol{N_i}$, then we can consider the restriction of $\rho$ to the segment joining 
$(x,y,t_1,t_2)$ with $q_i(s)$ and write its Taylor expansion of 
second order centred at $q_i(s)$. By (\ref{null}) and the fact that the
Hessian matrix of $\rho$ is negative definite at $q_i(0)$ (and then, by
continuity in a small neighbourhood), we have that there exist $\delta,C>0$
such that, if $U$ is small enough and $|t_1-u_{i-1}(0,0)|<\delta$,
$|t_2-u_{i}(0,0)|<\delta$, then
$$\rho(x,y,t_1,t_2)\leq 1-C(\nu_i\!\cdot\!(x,y))^2-C(t_1-u_{i-1}(s\tau_i))^2
-C(t_2-u_i(s\tau_i))^2\leq 1.$$
In the general case, since the curve $\{y=\psi_i(x,0) \}$ is $C^2$ with null
second derivative at $0$, one can find $s>0$, $a\in\real$ such that the
segment joining $(x,y)$ with $s\tau_i+a\nu_i$ is all contained in ${\ol N_i}$
and the ratio $|(x,y)-s\tau_i-a\nu_i|/a^2$ is infinitesimal as $a\to 0$.
Since $s>0$, the segment joining $s\tau_i+a\nu_i$ with its projection
$s\tau_i$ on $S_{i-1,i}$ is all contained in ${\ol N_i}$, so that we can apply
to this point the estimate above; if we call $L$ the $L^{\infty}$-norm of the
gradient of $\rho$, we obtain that
\begin{eqnarray*}
\rho(x,y,t_1,t_2) & \leq & \rho(s\tau_i+a\nu_i,t_1,t_2)
+L|(x,y)-s\tau_i-a\nu_i| \\
& \leq & 1- a^2\left(C -L\frac{|(x,y)-s\tau_i-a\nu_i|}{a^2}\right)
-C(t_1-u_{i-1}(s\tau_i))^2 -C(t_2-u_i(s\tau_i))^2,
\end{eqnarray*}
which is less than $1$, provided $U$ is small enough.
So we have proved that, if $\eps$ is sufficiently small, then there exists $\delta>0$ 
such that
\begin{equation}\label{ob1} 
\rho(x,y,t_1,t_2)\leq 1 \qquad \text{for}\ (x,y)\in \ol{N_i},  
\ |t_1-u_{i-1}(0,0)|<\delta, \ |t_2-u_i(0,0)|<\delta,
\end{equation}
provided $U$ is sufficiently small.

\vspace{.1cm}
Suppose now $(x,y)\in\ol{P_i}$, $|t_1-u_{i-1}(0,0)|\leq\eps$,
$|t_2-u_i(0,0)|\leq\eps$. In order to show that $\rho\leq 1$ also in this case,
we will compute the traces of the gradient and of the Hessian matrix of $\rho$
at the point $q_{i}(0)$. The main difference with respect to the previous case
is that in the region $\ol{P_i}$ the function $\beta_i-\alpha_i$ coincides with
the solution $\eta_i$ of the problem \begin{equation}\label{dfr2} \begin{cases}
\displaystyle
\frac{1}{\lambda}\sigma_i(x,y)\phi_i(x,y)
\cdot\nabla\eta_i(x,y) =\frac{\eps^2}{v_{i-1}^2(x,y)}
-\frac{\eps^2}{v_i^2(x,y)} , \\
\eta_i(s,0)= 0 \quad (s\geq 0),
\end{cases}
\end{equation}
while the function $\sigma_i$ is defined as
\begin{equation}\label{defs2}
\sigma_i(x,y)=
\frac{g(h_i(x,y))}{2\phi_i^y(h_i(x,y),0)} \qquad \forall
(x,y)\in \ol{P_i}.
\end{equation}
By (\ref{cpei}) and (\ref{mulgrad}) it follows that 
\begin{equation}\label{eq2}
\nabla \rho(q_{i}(0))= \nabla I^{\nu_i} (q_{i}(0)). 
\end{equation}
By (\ref{expr}) we obtain
the following expression for the gradient of $I^{\nu_i}$ with respect to
$\tau_i,\nu_i$ computed at the point $q_i(0)$: \begin{equation}\label{inter}
\nabla_{\!\tau_i\nu_i} I^{\nu_i}(q_i(0))=
g(0)\nabla\sigma_i(0,0)+\nabla\phi_i^{\nu_i}(0,0)+\frac{\sqrt{3}}{2}\tau_i,
\end{equation}
where we have used the Euler conditions, the fact that
$\nabla(\beta_i-\alpha_i)(0,0)=0$ by (\ref{dfr2}), and that
$$\nabla v_{i-1}(x,y)+\nabla v_i(x,y)=-\tau_i \qquad \forall (x,y)\in U.$$
It follows immediately by (\ref{prime}) that
\begin{equation}\label{336p}
\nabla\phi_i^{\nu_i}(x,y)=g'(\tau_i\!\cdot\!(x,y))\tau_i
\end{equation}
and by the
definition of $g$ that  
\begin{equation}\label{gder}
g'(t)  =  \sqrt{3}\eps^2\frac{\df_x v_0(t,0)}
{v^2_0(t,0)}  =  -\frac{\sqrt{3}}{2}\eps^2\frac{1}{v^2_0(t,0)}
\end{equation}
for all $t\in\real$.
By differentiating (\ref{defs2}), we obtain that
\begin{equation}\label{nblsig}
\nabla\sigma_i(x,y)=\frac{1}{2}p(h_i(x,y))\nabla h_i(x,y),
\end{equation}
where we have set
$$p(t):=\frac{g'(t)}{\phi_i^y(t,0)}-
\frac{g(t)}{[\phi_i^y(t,0)]^2}\df_x\phi_i^y(t,0).$$
To compute the gradient of $h_i$ it is enough to differentiate the
second equality in (\ref{hdef}): this provides
\begin{equation}\label{hder}
\partial_x\psi_i(x,h_i)+\partial_s\psi_i(x,h_i)\partial_x h_i=0, \qquad
\partial_s\psi_i(x,h_i)\partial_y h_i=1;
\end{equation}
by (\ref{psider}) we have that
\begin{equation}\label{339p}
\nabla h_i(0,0)=-2\tau_i. 
\end{equation}
Since 
$$\df_x\phi_i^y(x,y)=
(-1)^{i+1}\frac{3}{4}f'(\nu_i\!\cdot\!(x,y))-\frac{1}{4}g'(\tau_i\!\cdot\!(x,y)),$$ 
we find that $p(0)=3g'(0)/g(0)$, and substituting in (\ref{nblsig}), we have
that 
\begin{equation}\label{sigh}
\nabla\sigma_i(0,0) = -3\frac{g'(0)}{g(0)}\tau_i.
\end{equation}
Since the partial derivatives of $I^{\nu_i}$ with respect to $t_1$ and $t_2$
are still given by (\ref{26}), they are equal to $0$ at the point
$q_{i}(0)$, as in the previous case. 
Then, by (\ref{eq2}), (\ref{inter}), (\ref{336p}), (\ref{sigh}), and
(\ref{gder}), we deduce that  \begin{equation}\label{pll}
\nabla \rho(q_{i}(0))=
\left(\frac{3\sqrt{3}}{2}\tau_i,0,0\right).
\end{equation}
To conclude the study of $\rho$ in this region, we write the Hessian matrix of
$\rho$ with respect to  $\nu_{i},t_1,t_2$, which still satisfies
(\ref{grads}). 
Differentiating (\ref{exprt}) and using the Euler conditions, the fact that
$\nabla(\beta_i-\alpha_i)(0,0)=0$, $\phi_i^{\tau_i}(0,0)=0$ and (\ref{dntn}),
we obtain that (\ref{tg}) still holds.
Differentiating (\ref{fnd}) and computing the result at $q_i(0)$, we have that
\begin{equation}\label{341p}
\dfii I^{\nu_i}(q_i(0)) =
\frac{1}{\lambda}g(0)\dfii(\beta_i-\alpha_i)(0,0)+g(0)
\dfii\sigma_i(0,0)+\frac{3}{2\eps}(\dfi v_{i-1}(0,0)- 
\dfi v_i(0,0)),
\end{equation}
where we have used in particular that $\dfi \sigma_i(0,0)=0$ by (\ref{sigh})
and that $\dfii\phi_i^{\nu_i}\equiv 0$.
In order to compute the second derivative of $\beta_i-\alpha_i$ with respect to the direction $\nu_i$, 
we differentiate (\ref{dfr2}) with respect to $x$ and 
with respect to $y$; using the fact that $\df_x(\beta_i-\alpha_i)(s,0)=0$ for
every $s\geq 0$, we obtain 
\begin{eqnarray}\label{vana}
& \partial^2_{x}(\beta_i-\alpha_i)(0,0)=0,  \qquad
\label{mix}
\partial^2_{xy}(\beta_i-\alpha_i)(0,0)= \displaystyle \frac{6}{\eps}
(-1)^{i+1}\frac{\lambda}{g(0)}, & \vs \\
& \partial^2_{y}(\beta_i-\alpha_i)(0,0)  = \displaystyle
-\frac{2\sqrt{3}}{\eps}\frac{\lambda}{g(0)}+ \sqrt{3}(-1)^{i+1}
\partial^2_{xy}(\beta_i-\alpha_i)(0,0)
=  \frac{4\sqrt{3}}{\eps}\frac{\lambda}{g(0)}. & \label{pure}
\end{eqnarray}
By the relation
$\dfii=\frac{3}{4}\partial^2_x +\frac{\sqrt{3}}{2}(-1)^{i}\partial^2_{xy}
+\frac{1}{4}\partial^2_y$, it follows that
$$\dfii(\beta_i-\alpha_i)(0,0)=
-\frac{2\sqrt{3}}{\eps}\frac{\lambda}{g(0)}.$$
Since $\dfi h_i(0,0)=0$ by (\ref{339p}), from (\ref{nblsig}) we obtain that  
\begin{equation}\label{343p}
\dfii\sigma_i(0,0) =  \frac{1}{2}\left(
\frac{g'(0)}{\phi_i^y(0,0)}-\frac{g(0)}{[\phi_i^y(0,0)]^2}\df_x\phi_i^y(0,0)\right)\dfii h_i(0,0) 
=  \frac{3}{2}\frac{g'(0)}{g(0)}\dfii h_i(0).
\end{equation}
By differentiating twice with respect to the direction $\nu_i$ the second
equality in (\ref{hdef}), we obtain that
$$(\nu_i^x)^2 \partial^2_x\psi_i(x,h_i)+2\nu_i^x\partial^2_{xs}\psi_i(x,h_i)
\dfi h_i+ \partial^2_s\psi_i(x,h_i)(\dfi h_i)^2
+\partial_s\psi_i(x,h_i)\dfii h_i =0;$$
since $\dfi h_i(0,0)=0$ by (\ref{339p}) and $\df^2_x\psi_i(0,0)=0$ by
(\ref{psider2}), we can conclude that 
$\dfii h_i(0,0)= 0$ and then, by (\ref{343p}) also the limit of
$\dfii\sigma_i$ at $(0,0)$ is equal to $0$. Taking (\ref{veq}) and (\ref{341p})
into account,  we can conclude that  $$\dfii I^{\nu_i}(q_i(0))
=-\frac{\sqrt{3}}{2\eps},$$
i.e., (\ref{nm}) is still satisfied.
Since it is easy to see that also the other second derivatives of $\rho$
remain unchanged, we can conclude that the Hessian matrix of $\rho$ with
respect to $\nu_i,t_1,t_2$ is negative definite at $q_{i}(0)$.

If the segment joining $(x,y,t_1,t_2 )$
with $q_i(0)$ is all contained in $\ol{P_i}$, then we consider the Taylor
expansion of second order centred at $q_i(0)$ of the function $\rho$
restricted to this segment; since the component of $(x,y)$ along
$\tau_{i}$ is less or equal than $0$, by
(\ref{pll}) and by the fact that the Hessian matrix of $\rho$ with respect to
$\nu_{i},t_1,t_2$ is negative definite, we have that there exists $\delta>0$
such that $\rho(x,y,t_1,t_2)\leq 1$ for 
$|t_1-u_{i-1}(0,0)|<\delta$, $|t_2-u_{i}(0,0)|<\delta$, provided $U$ is small
enough.
In the general case, we can find $s\leq 0$, $a\in\real$ such that the segments
joining $(x,y)$ with $s\tau_i+a\nu_i$, and $s\tau_i+a\nu_i$ with $(0,0)$ are
all contained in $\ol{P_i}$, and $|(x,y)-s\tau_i-a\nu_i|/a^2$ is infinitesimal
as $a\to 0$. Arguing as for the region $\ol{N_i}$, this is enough to obtain the
same conclusion.
So we have proved that, if $\eps$ is small enough, there exists $\delta>0$
such that 
\begin{equation}\label{ob2}
\rho(x,y,t_1,t_2)\leq 1 \qquad \text{for}\ (x,y)\in 
\ol{P_i},   \ |t_1-u_{i-1}(0,0)|<\delta, \ |t_2-u_i(0,0)|<\delta,
\end{equation}
provided $U$ is sufficiently small.

By (\ref{ob1}) and (\ref{ob2}) Step~1 is proved. \qed

\section{Estimates for $t_1$ and $t_2$ near $u_0$ and $u_2$}

This section is devoted to the proof of the following step.

\bstep{2}{For a suitable choice of the function $f$ (see (\ref{phi})), there
exists $\delta >0$ such that condition (d) holds for
$|t_1-u_0(0,0)|<\delta$, $|t_2-u_2(0,0)|<\delta$, provided $U$ is small
enough.}

\noindent In order to prove the step, we want to show that the function $\rho$,
introduced at the beginning of Section~4, is less or equal than $1$ in a
neighbourhood of the point $(0,0,u_0(0,0),u_2(0,0))$.  
We can assume that $|t_1-u_0(0,0)|\leq\eps$, $|t_2-u_2(0,0)|\leq\eps$. Since
now the derivatives of $\rho$ may be discontinuous on the curves $\{h_1=0\}$
and $\{h_2=0\}$, we have to consider separately four different cases,
one for $(x,y)$ belonging to each one of the regions $N_1\cap N_2$, $N_1\cap
P_2$, $N_2\cap P_1$, and $P_1\cap P_2$.

Let $I^x$ and $I^y$ be the components of the integral in
(\ref{intgr}) with respect to $e^x$ and $e^y$, that are the
tangent and the normal direction, respectively, to the third part of the
discontinuity set $S_{0,2}$.

\vspace{.1cm}

Consider first the case $(x,y)\in \ol{P_1}\cap \ol{P_2}$, which is the
region containing $S_{0,2}$; as before, we will study the derivatives of
$\rho$ at the points of the form 
$$q_0(x):=(x,0,u_0(x,0),u_2(x,0)), \qquad x\geq 0.$$
Condition (\ref{e0}) implies that $\rho(q_0(x))=1$ for every $x\geq 0$; we
want to prove that 
\begin{equation}\label{nlo}
\nabla \rho(q_{0}(x))=0 \qquad \forall x\geq 0
\end{equation}
and that the Hessian matrix of $\rho$ with respect to $y,t_1,t_2$ is negative
definite at $q_i(0)$.
By the definition of $\rho$, it follows that
$$\nabla \rho=\frac{1}{\rho}(I^x \nabla I^x + I^y\nabla I^y).$$
Since $I^x(q_{0}(x))=0$ and $I^y(q_{0}(x))=1$ for every $x\geq 0$, we have that
$$\nabla \rho(q_{0}(x))=\nabla I^y(q_0(x)).$$
By (\ref{02}) and by the definition of $\vp$ in $G_i$ we can write the
explicit expression of $I^y$ at the point $(x,y,t_1,t_2)$: 
\begin{multline}\label{ycomp} 
I^y = -2(t_1-u_0)\dfo u_0+2(t_2-u_2)\dfo u_2+
\frac{1}{\lambda}\sum_{i=1}^2(\beta_i -\alpha_i
+\lambda)\sigma_i\phi_i^y \\
+\frac{\sqrt{3}}{2v_0}(\eps^2-(t_1-u_0)^2)+\frac{\sqrt{3}}{2v_2}(\eps^2-(t_2-u_2)^2),
\end{multline} 
and by differentiating with respect to $y$, we obtain 
\begin{multline}\label{fnd0}
\dfo I^y= 2(\dfo u_0)^2-2(\dfo u_2)^2
-2(t_1-u_0)\dfo^2 u_0+2(t_2-u_2)\dfo^2 u_2 \\
+\frac{1}{\lambda}\sum_{i=1}^2[\dfo(\beta_i-\alpha_i)\sigma_i \phi_i^y +
(\beta_i-\alpha_i+\lambda)\dfo(\sigma_i\phi_i^y)]  
-\frac{3}{4v_0^2}(\eps^2-(t_1-u_0)^2)\\
+\frac{3}{4v_2^2}(\eps^2-(t_2-u_2)^2) 
+\frac{\sqrt{3}}{v_0}(t_1-u_0)\dfo
u_0+\frac{\sqrt{3}}{v_2}(t_2-u_2)\dfo u_2.
\end{multline}
Since in the region $\ol{P_1}\cap \ol{P_2}$ the functions $\beta_i-\alpha_i$ 
coincide with the solutions of the problems (\ref{dfr2}), it results that
$\dfo(\beta_i-\alpha_i)(x,0)=0$ for $i=1,2$. Moreover, differentiating
(\ref{sigmai}) and the second equality in (\ref{phii})
with respect to $y$, we have that
\begin{equation}\label{rlt}
\dfo\sigma_2(x,y)=-\dfo\sigma_1(x,-y), \qquad
\dfo\phi_2^y(x,y)=-\dfo\phi_1^y(x,-y),
\end{equation}
and then, using again (\ref{phii}) and (\ref{sigmai}),
$$\phi_1^y(x,0)\dfo\sigma_1(x,0)=-\phi_2^y(x,0)\dfo\sigma_2(x,0), \quad
\sigma_1(x,0)\dfo\phi_1^y(x,0)=-\sigma_2(x,0)\dfo\phi_2^y(x,0).$$
By the Euler conditions, $\dfo u_0(x,0)=\dfo u_2(x,0)=0$ for every
$x\geq 0$; using all these remarks and (\ref{veq0}), we
deduce that  $\dfo I^{y}(q_0(x))=0$ for every $x>0$ and the 
equality holds also for the trace of $\dfo I^{y}$ at $q_0(0)$.
Since we have that
\begin{equation}\label{t1t2}
\df_{t_1} I^y=-2\dfo u_{0}-\frac{\sqrt{3}}{v_0}(t_1-u_0),\qquad \df_{t_2}
I^y= 2\dfo u_{2}-\frac{\sqrt{3}}{v_2}(t_2-u_2),
\end{equation}
by the Euler conditions it follows that $\df_{t_1} I^y(q_0(x))=\df_{t_2}
I^y(q_0(x))=0$. As $I^y(q_0(x))=1$ for every $x\geq 0$, this implies that
$\df_xI^y(q_0(x))=0$. Thus we have obtained equality (\ref{nlo}).

By (\ref{nlo}) and (\ref{e0})
the Hessian matrix of $\rho$ computed at $q_{0}(0)$ reduces to
\begin{equation}\label{zgrads}
\vhess \rho(q_{0}(0))=[\nablaz I^x\otimes \nablaz I^x+ \vhess
I^y](q_{0}(0)).
\end{equation}
As before, we know that
\begin{multline}\nonumber
I^x = -2(t_1-u_0)\df_x u_0+2(t_2-u_2)\df_x u_2+
\frac{1}{\lambda}\sum_{i=1}^2(\beta_i -\alpha_i
+\lambda)\sigma_i\phi_i^x \\
-\frac{1}{2v_0}(\eps^2-(t_1-u_0)^2)+\frac{1}{2v_2}(\eps^2-(t_2-u_2)^2),
\end{multline}
hence, by the Euler condition, the fact that
$\dfo(\beta_i-\alpha_i)(0,0)=0$ for $i=1,2$, and (\ref{nbl-0}), it
results that 
$$\dfo I^x(q_0(0)) = \frac{\sqrt{3}}{2} +
\sum_{i=1}^2\dfo(\sigma_i\phi_i^x)(0,0) =
\frac{\sqrt{3}}{2} +2\dfo\phi_1^x(0,0) + 2
\phi_1^x(0,0)\dfo\sigma_1(0,0),$$
where we have also used the first equalities in (\ref{phii}) and in
(\ref{rlt}), and the relation $\dfo\phi_2^x(x,y)=\dfo\phi_1^x(x,-y)$.
From (\ref{sigh}) we obtain that
$$\dfo\sigma_1(0,0)=\frac{3\sqrt{3}}{2}\frac{g'(0)}{g(0)}.$$
Then, using the definition of $\phi_1^x$ and (\ref{gder}), we can conclude that
\begin{equation}\label{tg0}
\dfo I^x(0,0)=\frac{\sqrt{3}}{2} -3g'(0)=
2\sqrt{3}.
\end{equation}
By differentiating (\ref{fnd0}) with respect to $y$ and by using the Euler
condition and the fact that $\dfo(\beta_i-\alpha_i)(0,0)=0$ for $i=1,2$, we obtain
$$\df^2_y I^y(q_0(0))= \frac{1}{\lambda}\sum_{i=1}^2[\df^2_y(\beta_i-\alpha_i)\phi_i^y
+\df^2_y(\sigma_i\phi_i^y)](0,0)+\frac{3\sqrt{3}}{2\eps}.$$
Equality (\ref{pure}) implies that
\begin{equation}\label{noph}
\frac{1}{\lambda}\sum_{i=1}^2[\df^2_y(\beta_i-\alpha_i)\sigma_i\phi_i^y](0,0)=
\frac{4\sqrt{3}}{\eps}.
\end{equation}
In order to write explicitly $\df^2_y\sigma_i$ at $(0,0)$, we
differentiate the $y$-component in (\ref{nblsig}) with respect to $y$ and
we pass to the limit, taking into account that $\dfo
h_i(0)=(-1)^{i+1}\sqrt{3}$ by (\ref{339p}): 
$$\df^2_y\sigma_1(0,0) =\frac{3}{2}p'(0) +\frac{1}{2}p(0)\dfo^2 h_i(0).$$ 
By differentiating with respect to $y$ the second equality
in (\ref{hder}), we obtain that
$$\df^2_y h_1(0,0)=-(\dfo
h_1(0,0))^2\frac{\df^2_s\psi_1(0,0)}{\df_s\psi_1(0,0)}=0,$$ 
where the last equality follows by (\ref{psider2bis}).
Since
\begin{equation}\label{ppr}
p'(0)=2\frac{g''(0)}{g(0)}+3 \frac{[g'(0)]^2}{g^2(0)}
-4\frac{\df^2_x \phi_1^y(0,0)}{g(0)},
\end{equation}
while
\begin{equation}\label{foll}
\df^2_x\phi_1^y(0,0)=-\frac{3\sqrt{3}}{8}f''(0)+\frac{1}{8}g''(0),
\quad \df^2_y\phi_1^y(0,0)=-\frac{\sqrt{3}}{8}f''(0)+
\frac{3}{8}g''(0),
\end{equation}
and $g''(0)= -\sqrt{3}/(2\eps)$, we can write that
\begin{eqnarray}
\frac{1}{\lambda}\sum_{i=1}^2(\beta_i-\alpha_i+\lambda)
\displaystyle \df^2_y(\sigma_i\phi_i^y)(0,0) & = & 
(2\phi_1^y\df^2_y\sigma_1+
4\dfo \sigma_1\dfo \phi_1^y
+2\df^2_y\phi_1^y)(0,0) \nonumber \\
& = & 2\sqrt{3}f''(0)+3g''(0) \nonumber \\
& = & 2\sqrt{3}f''(0)-\frac{3\sqrt{3}}{2\eps}.\label{qnoph}
\end{eqnarray}
Substituting (\ref{noph}) and (\ref{qnoph}) in the expression of $\dfo^2 I^y$,
we find that 
\begin{equation}\label{znm}
\df^2_y I^y(q_0(0))= 2\sqrt{3}f''(0)+\frac{4\sqrt{3}}{\eps}.
\end{equation}
From (\ref{zgrads}), (\ref{tg0}), and (\ref{znm}), we finally obtain that
\begin{equation}\label{final0}
\df^2_y \rho(q_{0}(0))=[\dfo I^x(q_0(0))]^2+\dfo^2 I^y(q_0(0))=
12  +\frac{4\sqrt{3}}{\eps}+
2\sqrt{3}f''(0).
\end{equation}
As in the previous step, we can compute explicitly the other 
elements of the Hessian matrix of $\rho$ and we find that
$$\det\left(\begin{array}{cc} \dfo^2\rho & \df^2_{yt_1}\rho \vspace{.1cm}\\
\df_{yt_1}\rho & \df^2_{t_1}\rho \end{array} \right)(q_0(0)) =
-\frac{6}{\eps}f''(0)
-\frac{12\sqrt{3}}{\eps} -\frac{12}{\eps^2} -
4(\dfo^2u_0(0,0))^2,$$
$$\det\vhess \rho(q_0(0)) = \frac{6\sqrt{3}}{\eps^2}f''(0)+
\frac{36}{\eps^2}+\frac{12\sqrt{3}}{\eps^3}+
\frac{4\sqrt{3}}{\eps}[(\dfo^2u_0(0,0))^2+(\dfo^2u_2(0,0))^2].$$
If we impose the following condition on the second derivative of $f$ at $0$:
\begin{equation}\label{condf}
f''(0)<-2\sqrt{3}-\frac{2}{\eps}-\frac{2\eps}{3}
[(\dfo^2u_0(0,0))^2+(\dfo^2u_2(0,0))^2],
\end{equation}
then the Hessian matrix of $\rho$ is negative 
definite at $q_{0}(0)$. 

To conclude, we restrict $\rho$ to the segment joining $(x,y,t_1,t_2)$ with
$q_0(x)$ and we write its Taylor expansion of second order centred at
$q_0(x)$; using (\ref{nlo}) and choosing $f$ satisfying (\ref{condf}) (so that
the Hessian matrix of $\rho$ is negative  definite at $q_{0}(0)$, and then by
continuity in a small neighbourhood), we obtain that there exists $\delta>0$
such that \begin{equation}\label{ob01} \rho(x,y,t_1,t_2)\leq 1 \qquad
\text{for}\ (x,y)\in  \ol{P_1}\cap \ol{P_2},  \ |t_1-u_{0}(0,0)|<\delta, \
|t_2-u_2(0,0)|<\delta, 
\end{equation}
provided $U$ is sufficiently small.

\vspace{.1cm}

Let us consider the set $\ol{N_1}\cap \ol{N_2}$: in this region 
$\sigma_1=\sigma_2=1$, while the functions $\beta_i-\alpha_i$ 
coincide with the solutions of the problems (\ref{dfr}). 
By (\ref{e0}) the gradient of $\rho$ at the point $q_0(0)$ is given by
\begin{equation}\label{eq02}
\nabla \rho(q_0(0))=\nabla I^y(q_0(0)).
\end{equation}
By (\ref{ycomp}) we derive the
explicit expression for the gradient of $I^y$ with respect to $x,y$;
using the Euler condition, the fact that $\nabla(\beta_i-\alpha_i)(0,0)=0$,
the constancy of $\sigma_i$ and the equality
\begin{equation}\label{360p}
\nabla v_0(x,y)+\nabla v_2(x,y)=-\tzer \qquad \forall (x,y)\in U,
\end{equation}
we obtain that
$$\nabla_{\! xy}I^y(q_0(0))= \sum_{i=1}^2\nabla \phi_i^y(0,0)
+ \frac{\sqrt{3}}{2}\tzer
=-\frac{1}{2}g'(0)\tzer+\frac{\sqrt{3}}{2}\tzer=
\frac{3\sqrt{3}}{4}\tzer.$$
Since the partial
derivatives of $I^y$ with respect to $t_1$ and $t_2$ are still given by
(\ref{t1t2}), they are equal to $0$ at $q_0(0)$, as in the previous case.
Therefore, we have that
\begin{equation}\label{posi}
\nabla \rho(q_0(0))=\left(\frac{3\sqrt{3}}{4} \tzer,0,0\right).
\end{equation}

If $(x,y)\neq (0,0)$ belongs to $\ol{N_1}\cap \ol{N_2}$ and
the segment joining $(x,y)$ with $(0,0)$ is all contained in
$\ol{N_1}\cap \ol{N_2}$, then by the Mean Value Theorem,
(\ref{posi}) and the fact that $x$ is strictly negative, we can conclude
that there exists $\delta>0$ such that
\begin{equation}\label{ob02}
\rho(x,y,t_1,t_2)\leq 1 \qquad \text{for}\ |t_1-u_{0}(0,0)|<\delta, \
|t_2-u_2(0,0)|<\delta, 
\end{equation}
provided $U$ is sufficiently small.
If the segment joining $(x,y)$ with $(0,0)$ is not contained in $\ol{N_1}\cap
\ol{N_2}$, then we can find a regular curve connecting $(x,y)$ and
$(0,0)$, along which we can repeat the same estimate as above.

\vspace{.1cm}

At last consider the set $\ol{N_2}\cap \ol{P_1}$, since the case
$\ol{N_1}\cap \ol{P_2}$ is completely analogous. 
In this region, $\sigma_1$ is defined by (\ref{defs2}), while
$\sigma_2$ is identically equal to $1$; the function $\beta_1-\alpha_1$
coincides with the solution of the problem (\ref{dfr2}) for $i=1$, while
$\beta_2-\alpha_2$ with the one of (\ref{dfr}) for $i=2$. 
Equality (\ref{eq02}) still holds, as well as the fact
that $\nabla(\beta_i-\alpha_i)(0,0)=(0,0)$ for all $i$; since
$\nabla\sigma_1$ is given by the formula (\ref{sigh}) and
$\nabla\sigma_2\equiv 0$, by (\ref{phi}), (\ref{gder}), (\ref{ycomp}), and
(\ref{360p}) we have that  \begin{eqnarray*}
\nabla_{\! xy} I^y(q_0(0))& = & \sum_{i=1}^2\nabla \phi_i^y(0,0)
+\phi_1^y(0,0)\nabla \sigma_1(0,0)+\frac{\sqrt{3}}{2}\tzer \\
& = & \frac{3\sqrt{3}}{4}(\tzer+\tone)= -\frac{3\sqrt{3}}{4}\ttwo,
\end{eqnarray*}
hence
$$\nabla \rho(q_0(0))=\left(-\frac{3\sqrt{3}}{4} \ttwo,0,0 \right).$$ 
Since the gradient of $\rho$ vanishes along the direction
$(\ntwo,0,0)$, we need to compute the Hessian matrix of $\rho$ with respect to
$\ntwo,t_1,t_2$ at the point $q_0(0)$; from the equality
$\nabla_{\!\ntwo t_1 t_2}I^y(q_0(0))=0$, we have that
\begin{equation}\label{drag}
\vvhess \rho(q_0(0))=[\nabla_{\!\ntwo t_1t_2}I^x \otimes \nabla_{\!\ntwo
t_1t_2}I^x + \vvhess I^y](q_0(0)).
\end{equation}
Using the fact that $\nabla u_0(0,0)=\nabla u_2(0,0)=0$ and
$\nabla(\beta_i-\alpha_i)(0,0)=0$, we obtain
\begin{eqnarray*}
\df_{\nu_2}I^x(q_0(0)) & = & \sum_{i=1}^2 \df_{\nu_2}\phi_i^x(0,0)+
\df_{\nu_2}\sigma_1(0,0)\phi_1^x(0,0)+\frac{1}{2} \df_{\nu_2}(v_0-v_2) \\
& = & \dfo\phi_1^x(0,0)-\frac{9}{4}g'(0)+\frac{\sqrt{3}}{4}=
\sqrt{3},
\end{eqnarray*}
where the second equality follows from (\ref{sigh}) and from the fact that 
$\df_{\nu_2}\phi_1^x+\df_{\nu_2}\phi_2^x=\dfo\phi_1^x$ at $(0,0)$.
If we differentiate (\ref{ycomp}) twice with respect to the direction $\ntwo$
and we compute the result at the point $q_0(0)$, we obtain 
\begin{equation}\label{bigder}
\df^2_{\nu_2}I^y(0,0)=\left(\frac{1}{\lambda}\sum_{i=1}^2\df^2
_{\nu_2}(\beta_i-\alpha_i)\sigma_i\phi_i^y +\sum_{i=1}^2\df^2_{\nu_2}\phi_i^y
+\df^2_{\nu_2}\sigma_1\phi_1^y+2\df_{\nu_2}\sigma_1\df_{\nu_2}\phi_1^y\right)(0,0) +
\frac{3\sqrt{3}}{4\eps}.
\end{equation}
From (\ref{mix}) and (\ref{pure}), and from (\ref{star}) it follows
respectively that
\begin{equation}\label{noteq}
\df^2_{\nu_2}(\beta_1-\alpha_1)(0,0)=
\frac{4\sqrt{3}}{\eps^3}\frac{\lambda}{g(0)},\quad
\df^2_{\nu_2}(\beta_2-\alpha_2)(0,0)=
-\frac{2\sqrt{3}}{\eps}\frac{\lambda}{g(0)}.
\end{equation}
Since by (\ref{nblsig}) we have that $\df_{\ntwo}\sigma_1(x,y)=
\frac{1}{2}p(h_1(x,y))\df_{\ntwo}h_1(x,y)$, then 
$$\df^2_{\ntwo}\sigma_1(0,0)= \frac{1}{2}p'(0)(\df_{\ntwo}h_1(0,0))^2+
\frac{1}{2}p(0)\df^2_{\ntwo}h_1(0,0).$$
Some easy computations show that $\df^2_{\ntwo}h_1(0,0)=0$; using (\ref{339p})
it results that \begin{equation}\label{2ds}
\df^2_{\ntwo}\sigma_1(0,0)= \frac{3}{2}p'(0)=
\frac{9}{2}\frac{[g'(0)]^2}{g^2(0)}+\frac{9}{4}\sqrt{3}\frac{f''(0)}{g(0)},
\end{equation}
where the last equality follows by (\ref{ppr}) and by the first equality in
(\ref{foll}).
At last, by using (\ref{phii}) and (\ref{foll}), we obtain that
\begin{equation}\label{last}
\sum_{i=1}^2\df^2_{\nu_2}\phi_i^y(0,0)=
\frac{3}{4}\df^2_x\phi_1^y(0,0)+\frac{1}{4}\df^2_y\phi_1^y(0,0)=
-\frac{5}{8}\sqrt{3}f''(0)+\frac{3}{8}g''(0),
\end{equation}
and by substituting (\ref{noteq}), (\ref{2ds}), and (\ref{last}) in
(\ref{bigder}), we deduce that
$$\df^2_{\nu_2}I^y(q_0(0))=
\frac{\sqrt{3}}{2}f''(0)+\frac{\sqrt{3}}{\eps},$$
hence
$$\df^2_{\ntwo}\rho(q_0(0))=
3+\frac{\sqrt{3}}{\eps}+\frac{\sqrt{3}}{2}f''(0).$$
By differentiating (\ref{t1t2}) with respect to $\ntwo$ and by (\ref{drag}), we
obtain 
$$\df^2_{\ntwo t_1}\rho(q_0(0))=-2\df_{\ntwo}\dfo u_0(0,0)=
-\dfo^2 u_0(0,0), \qquad \df^2_{\ntwo t_2}\rho(q_0(0))=2\df_{\ntwo}\dfo
u_2(0,0)= \dfo^2 u_2(0,0).$$
At this point, it is easy to see that, if $f$ satisfies the condition
\begin{equation}\label{condf2}
f''(0)<-2\sqrt{3}-\frac{2}{\eps}-\frac{\eps}{6}[(\dfo^2u_0(0,0))^2
+(\dfo^2u_2(0,0))^2]
\end{equation}
then the Hessian matrix of $\rho$ with respect to $\ntwo, t_1,t_2$ is negative
definite at the point $q_0(0)$. 
Arguing as for the region $\ol{P_i}$ in the previous section, 
it can be proved that, if $f$ satisfies (\ref{condf2}), then there
exists $\delta>0$ such that 
\begin{equation}\label{ob03}
\rho(x,y,t_1,t_2)\leq 1 \qquad \text{for}\ (x,y)\in
\ol{N_2}\cap \ol{P_1},  \ |t_1-u_{0}(0,0)|<\delta, \
|t_2-u_2(0,0)|<\delta,
\end{equation}
provided $U$ is sufficiently small.

\vspace{.1cm}

Since condition (\ref{condf}) implies (\ref{condf2}), if we require that
(\ref{condf}) holds, then
by (\ref{ob01}),
(\ref{ob02}), and (\ref{ob03}), we can conclude that Step~2 is true. \qed

\section{Proof of condition (d)}

In this section we complete the proof of condition (d). To this aim 
it is enough to check condition (d) in the three cases studied in the
following step, as it will be clear at the end of the section.

\bstep{3}{If $\eps$ is sufficiently small, $\delta\in (0,\eps)$, and $U$ 
is sufficiently small, condition (d) is true for $t_1\leq t_2$ whenever one
of the following three conditions is satisfied: 
\begin{description}
	\item[1)] $|t_1-u_0(0,0)|\geq\delta$ and $|t_1-u_1(0,0)|\geq\delta$;
	\item[2)] $|t_2-u_1(0,0)|\geq\delta$ and $|t_2-u_2(0,0)|\geq\delta$;
	\item[3)] $|t_1-u_0(0,0)|\geq\delta$ and $|t_2-u_2(0,0)|\geq\delta$.
\end{description}}

\noindent
Let us fix $\delta\in (0,\eps)$ and set
\begin{multline}\nonumber
M_1(x,y):=\max\{ |I(x,y,t_1,t_2)|: \ u_0(x,y)-\eps\leq t_1\leq t_2\leq
u_2(x,y)+\eps,\\
|t_1-u_0(0,0)|\geq\delta,\ \ |t_1-u_1(0,0)|\geq\delta \}.
\end{multline}
It is easy to see that the function $M_1$ is continuous. Let us prove that
$M_1(0,0)<1$. For simplicity of notation, from now on we will denote 
$I(0,0,t_1,t_2)$ simply by $I(t_1,t_2)$ and $u_i(0,0)$ by $u_i$.

Let $t_1,t_2$ be such that $u_0-\eps\leq t_1\leq t_2\leq u_2+\eps$ with
$|t_1-u_0|\geq\delta$ and $|t_1-u_1|\geq\delta$. Suppose furthermore that
$|t_1-u_1|\leq\eps$; then, we can write 
\begin{eqnarray*}
I(t_1,t_2) & = & I(t_1,u_1)+I(u_1,u_2)+I(u_2,t_2), \\
I(u_2,t_2) & = & I(u_2,t_2\lor(u_2-\eps))+I(u_2-\eps,t_2\land(u_2-\eps)).
\end{eqnarray*}
Therefore, we have
\begin{equation}\label{decomp}
I(t_1,t_2) =I(t_1,u_1)+I(u_1,u_2)+I(u_2,t_2\lor(u_2-\eps))
-I(t_2\land(u_2-\eps),u_2-\eps).
\end{equation}
From the definition of $\vp$ in $G_1,G_2$ it follows that 
\begin{eqnarray}\nonumber
& I(s_1,u_1)=\displaystyle -\frac{1}{\eps}(s_1-u_1)^2\tzer & \text{for} \
|s_1-u_1|\leq\eps, \\ & I(u_2,s_2)=\displaystyle 
\frac{1}{\eps}(s_2-u_2)^2\tone & \text{for} \ |s_2-u_2|\leq\eps;
\label{ovvio2}  \end{eqnarray}
using condition (e), we have that
\begin{equation}\label{rombino}
I(t_1,u_1)+I(u_1,u_2)+I(u_2,t_2\lor(u_2-\eps))\in
\ntwo-\frac{\delta^2}{\eps}\tzer+R_1,
\end{equation}
where $R_1$ is the parallelogram spanned by the vectors
$\eps\tone$ and $-\left(\eps-\frac{\delta^2}{\eps}
\right)\tzer$.
Let $C$ be the intersection of the half-plane
$\{(x,y)\in\real^2: \ntwo\!\cdot\!(x,y)\geq 1-\sqrt{3}\eps\}$ with the open
ball centred at $0$ with radius $1$; some elementary geometric considerations
show that \begin{equation}\label{half}
\ntwo-\frac{\delta^2}{\eps}\tzer+R_1 \, \subset \, C.
\end{equation}
If $T_i$ is the segment joining $0$ with $g(0)\nu_i$, then from the definition 
of $\vp$ in $K_i$, it follows that
\begin{equation}\label{bu}
I(u_{i-1}+\eps,u_i-\eps)=g(0)\nu_i,
\end{equation}
and
\begin{equation}\label{beau}
I(s_1,s_2)\in T_i
\end{equation}
for $u_{i-1}+\eps \leq s_1\leq s_2\leq u_i-\eps$, $i=1,2$.
Let $D:=-T_2$; from (\ref{decomp}), (\ref{rombino}), (\ref{half}),
and (\ref{beau}), we deduce that $$I(t_1,t_2) \in C+D;$$
since $g(0)=1-\sqrt{3}\eps$, the set $C+D$ is contained in the open
ball centred at $0$ with radius $1$. This concludes the proof when
$|t_1-u_1|\leq\eps$.

If $|t_2-u_1|\leq\eps$, we consider the decomposition
\begin{eqnarray*}
I(t_1,t_2) & = & I(t_1,u_0)+I(u_0,u_1)+I(u_1,t_2), \\
I(t_1,u_0) & = & I(t_1\land(u_0+\eps),u_0)+I(t_1\lor(u_0+\eps),u_0+\eps),
\end{eqnarray*}
and the proof is completely analogous.
 
When $|t_1-u_1|>\eps$ and $|t_2-u_1|>\eps$, we can write
\begin{eqnarray*}
I(t_1,t_2) & = & I(t_1,u_0)+I(u_0,u_2)+I(u_2,t_2), \\
I(t_1,u_0) & = & I(t_1\land(u_0+\eps),u_0)+I(t_1\lor(u_0+\eps),u_0+\eps), \\
I(u_2,t_2) & = & I(u_2,t_2\lor(u_2-\eps))+I(u_2-\eps,t_2\land(u_2-\eps));
\end{eqnarray*}
therefore, we have
\begin{multline}\label{decomp2}
I(t_1,t_2)=I(t_1\land(u_0+\eps),u_0)+I(u_0,u_2)+I(u_2,t_2\lor(u_2-\eps)) \\
+I(t_1\lor(u_0+\eps),t_2\land(u_2-\eps)) - I(u_0+\eps,u_2-\eps).
\end{multline}
Since from the definition of $\vp$ in $G_0$ it follows that 
\begin{equation}\label{ovvio0}
I(s_0,u_0)=-\frac{1}{\eps}(s_0-u_0)^2\ttwo \qquad \text{for} \
|s_0-u_0|\leq\eps,  \end{equation}
using condition (e) and (\ref{ovvio2}), we have that
\begin{equation}\label{rombetto}
I(t_1\land(u_0+\eps),u_0)+I(u_0,u_2)+I(u_2,t_2\lor(u_2-\eps))\in
\nzer-\frac{\delta^2}{\eps}\ttwo +R_2, 
\end{equation}
where $R_2$ is the parallelogram spanned by the
vectors $\eps\tone$ and $-\left(\eps-\frac{\delta^2}{\eps}
\right)\ttwo$.
Let $E$ be the parallelogram having as consecutive sides $T_1$ and $T_2$, and
let $F$ be the set $E-g(0)\nzer$; as $I(u_1-\eps,u_1+\eps)=0$, from (\ref{bu})
it follows that \begin{equation}\label{bel}
I(u_0+\eps,u_2-\eps)=g(0)\nzer=(1-\sqrt{3}\eps)\nzer,
\end{equation}
and from (\ref{beau}),
\begin{equation}\label{-bel}
I(s_1,s_2)\in E
\end{equation}
for every $u_0+\eps\leq s_1\leq s_2\leq u_2-\eps$, with $|s_1-u_1|>\eps$ and
$|s_2-u_1|>\eps$. 
From (\ref{decomp2}), (\ref{rombetto}), (\ref{bel}), (\ref{-bel}), we obtain
that
$$I(t_1,t_2)\in \nzer-\frac{\delta^2}{\eps}\ttwo +R_2+F.$$
The set $\nzer-\frac{\delta^2}{\eps}\ttwo +R_2+F$ is a polygon, since it is
the sum of two polygons, and it is possible to prove that, if 
$\eps<\sqrt{3}$, its vertices are all
contained in the open ball with centre $0$ and radius $1$. Then,
under this condition, the whole set $\nzer-\frac{\delta^2}{\eps}\ttwo
+R_2+F$ is contained in this ball; this concludes the proof of the
inequality $M_1(0,0)<1$. 

By continuity, choosing $U$ small enough, we obtain that $M_1(x,y)<1$ for every
$(x,y)\in U$, which proves 1).

To prove 2) and 3), we define analogously
\begin{multline}\nonumber
M_2(x,y):=\max\{ |I(x,y,t_1,t_2)|: \ u_0(x,y)-\eps\leq t_1\leq t_2\leq
u_2(x,y)+\eps,  \\
|t_2-u_1(0,0)|\geq\delta, \ \ |t_2-u_2(0,0)|\geq\delta \},
\end{multline}
\begin{multline}\nonumber
M_3(x,y):=\max\{ |I(x,y,t_1,t_2)|: \ u_0(x,y)-\eps\leq t_1\leq t_2\leq
u_2(x,y)+\eps, \\ 
|t_1-u_0(0,0)|\geq\delta, \ \ |t_2-u_2(0,0)|\geq\delta \}.
\end{multline}
It is easy to see that the functions $M_2$ and
$M_3$ are continuous and, arguing as in the case of $M_1$, we can prove that
$M_2(0,0)<1$ and $M_3(0,0)<1$, which yield 2) and 3) by continuity. \qed

\vspace{.5cm}\noindent

{\sc Conclusion.-- } As in Step~3, we simply write $u_i$ instead
of $u_i(0,0)$. Let us show that, if $f$ satisfies (\ref{condf}), and $\eps$
and $U$ are sufficiently small, then condition (d) is true for $u_0(x,y)-\eps
\leq t_1 < t_2\leq u_2(x,y)+\eps$ and in fact for every
$t_1,t_2\in\real$, since $\vpxi(x,y,z)=0$ for $z\leq u_0(x,y)-\eps$ and for
$z\geq u_2(x,y)+\eps$.

We start by considering the case $|t_1-u_0|<\delta$. If $|t_2-u_1|<\delta$,
the conclusion follows from Step~1. If $|t_2-u_1|\geq\delta$, the result is a
consequence of Step~2 when $|t_2-u_2|<\delta$, and of Step~3.2) in the other
case.

We consider now the case $|t_1-u_0|\geq\delta$. If $|t_1-u_1|\geq\delta$, the
conclusion follows from Step~3.1). If $|t_1-u_1|<\delta$, the result is a
consequence of Step~1 when $|t_2-u_2|<\delta$, and of Step~3.3) in the other
case.

\vspace{.2cm}
This concludes the proof of condition (d) and then, of Theorem~\ref{tj} in the
case $u_0$ symmetric. \qed

\section{The antisymmetric case}

In this section we show how the construction of the calibration for $u_i$
symmetric can be adapted to the antisymmetric case. 

If the function $u_0$ is antisymmetric with respect
to the bisecting line of $A_0$, then the reflection of $u_0$ with respect
to the $S_{0,1}$ and to $S_{0,2}$ provides an extension of $u_0$, which is
harmonic only on $\Omega\setminus S_{1,2}$ and which is multi-valued on
$S_{1,2}$, since the traces of the tangential derivatives of $u_0$ on
$S_{1,2}$ have different signs.  Since $u_1,u_2$ coincide, up to the sign and
to additive constants, with the reflections of $u_0$ with respect to $S_{0,1}$
and $S_{0,2}$, respectively, they are antisymmetric with respect to the
bisecting line of $A_1$ and $A_2$, respectively, and then, their extensions by
reflection are harmonic only on $\Omega\setminus S_{0,2}$ and
$\Omega\setminus S_{0,1}$, respectively. 

The calibration $\varphi$ can be defined as before, just
replacing the sets $G_0, G_1,G_2$ with
\begin{eqnarray*}
\tilde{G}_0 & = & \{(x,y,z)\in (U\setminus  S_{1,2}){\times}\real:
u_0(x,y)-\eps<z<u_0(x,y)+\eps \}, \\
\tilde{G}_1 & = & \{(x,y,z)\in (U\setminus  S_{0,2}){\times}\real:
u_1(x,y)-\eps<z<u_1(x,y)+\eps \}, \\
\tilde{G}_2 & = & \{(x,y,z)\in (U\setminus  S_{0,1}){\times}\real:
u_2(x,y)-\eps<z<u_2(x,y)+\eps \}, 
\end{eqnarray*}
and the sets $H_1,H_2$ with
\begin{eqnarray*}
\tilde{H}_1 & = & \{(x,y,z)\in (U\setminus  (S_{1,2}\cup
S_{0,2})){\times}\real: l_1+\lambda/2<z<l_1+3\lambda/2 \}, \\
\tilde{H}_2 & = & \{(x,y,z)\in (U\setminus (S_{0,1}\cup S_{0,2})){\times}\real:
l_2+\lambda/2<z<l_2+3\lambda/2 \}.
\end{eqnarray*}
Since $u_0$ is harmonic in $\Omega\setminus S_{1,2}$, the field $\varphi$ is
divergence free in $\tilde{G}_0$ by Lemma \ref{folio}.
Moreover, the normal component of $\varphi$ is continuous across the boundary
of $G_0$ since $\df_{\nu_2}u_0=\df_{\nu_2}v_0=0$ on $S_{1,2}$. The same
argument works for the sets $\tilde{G}_1,\tilde{G}_2$. 
By the harmonicity of $u_0$ and $u_1$, the field is divergence free in
$\tilde{H}_1$ and the normal component of $\varphi$ is continuous across the
boundary of $H_1$ since $\df_{\nu_2}u_0=0$ on $S_{1,2}$ and $\df_{y}u_1=0$ on
$S_{0,2}$.
Therefore, condition (a) is still satisfied in the sense of distributions on
$U{\times}\real$.

It is easy to see that conditions (b), (c), and (e) are satisfied.

The proof of Step~1, Step~2, and Step~3 can be easily adapted; indeed, even if
now the function $|I(x,y,t_1,t_2)|$ may present some discontinuities when
$(x,y)\in S_{i,j}$, we can write $U$ as the union of finitely many
Lipschitz open subsets $U_i$ such that $|I|$ is $C^2(\ol{U_i}{\times}\real^2)$
and study the behaviour of $|I|$ separately in each $\ol{U_i}$. So, it results
that also condition (d) is true. \qed

\end{document}